\documentclass{amsart}
\usepackage{color}
 
\usepackage{amsmath} 
\usepackage{amssymb}
\usepackage{mathrsfs}

\usepackage[dvipsnames]{xcolor}  

\newtheorem{theorem}{Theorem}[section]
\newtheorem{claim}[theorem]{Claim}

\newtheorem{lemma}[theorem]{Lemma} 
\newtheorem{claim/definition}[theorem]{Claim/Definition} 
 
\newtheorem{observation}[theorem]{Observation}

\theoremstyle{definition}
\newtheorem{definition}[theorem]{Definition}

\newtheorem{convention}[theorem]{Convention}
\newtheorem{fact}[theorem]{Fact}

\newtheorem{problem}[theorem]{Problem}
\newtheorem{discussion}[theorem]{Discussion}

\theoremstyle{remark}
\newtheorem{remark}[theorem]{Remark}

\newtheorem{notation}[theorem]{Notation}

\newcommand{\then}{{\underline{then}}}

\newcommand{\when}{{\underline{when}}}

\newcommand{\Then}{{\underline{Then}}}

\newcommand{\Iff}{{\underline{iff}}}

\newcommand{\wilog}{{\rm without loss of generality}}
\newcommand{\Wilog}{{\rm Without loss of generality}}

\newcommand{\rest}{{\restriction}}
\newcommand{\dom}{{\rm dom}} 
\newcommand{\true}{{\rm true}} 
\newcommand{\false}{{\rm false}} 
\newcommand{\cov}{{\rm cov}}
\newcommand{\meagre}{{\rm meagre}}

\newcommand{\crit}{{\rm crit}}

\newcommand{\GCH}{{\rm GCH}}

\newcommand{\bd}{{\rm bd}}
\newcommand{\otp}{{\rm otp}}

\newcommand{\Ord}{{\rm Ord}}

\newcommand{\lqq}{{``}}
 
\newcommand{\exi}{{\kappa }}
\newcommand{\exnu}{{ (\bar{ \nu } , \mathbf{G} (+)) }}

\newcommand{\dc}{{(< \lambda,\theta _ *, \bar{\theta } ) }}
\newcommand{\zdc}{{\boxdot _{< \lambda,\theta _ *, \bar{\theta } } }}  
\newcommand{\dic}{{directed complete }} 

\newcommand{\pigyon }{{\dagger  }}
\newcommand{\ggk}{{ \kappa   }} 
\newcommand{\ttg}{{\eta '   }}  
\newcommand{\ggy}{{ \eta   }}

\newcommand{\bfG}{{\mathbf G}}

\newcommand{\bfm}{{\mathbf m}}

\newcommand{\bfj}{{\mathbf j}}

\newcommand{\bfV}{{\mathbf V}}

\newcommand{\bfM}{{\mathbf M}}
\newcommand{\bfB}{{\mathbf B}}

\newcommand{\bfn}{{\mathbf n}}

\newcommand{\bfq}{{\mathbf q}}
\newcommand{\bfQ}{{\mathbf Q}}

\newcommand{\bbA}{{\mathbb A}}

\newcommand{\gb}{{\mathfrak b}}

\newcommand{\cU}{{\mathscr U}}

\newcommand{\bbD}{{\mathbb D}}

\newcommand{\bbR}{{\mathbb R}}

\newcommand{\cH}{{\mathscr H}}

\newcommand{\cF}{{\mathscr F}}

\newcommand{\gd}{{\mathfrak d}}

\newcommand{\bbL}{{\mathbb L}}

\newcommand{\gc}{{\mathfrak c}}

\newcommand{\bbP}{{\mathbb P}}
\newcommand{\cP}{{\mathscr P}}

\newcommand{\bbQ}{{\mathbb Q}}

\newcommand{\cS}{{\mathscr S}}

\newcommand{\varp}{{\varepsilon}}

\newcommand{\cX}{{\mathscr X}}
\newcommand{\cW}{{\mathscr W}}

\newcommand{\cf}{{\rm cf}}

\newcount\skewfactor
\def\mathunderaccent#1#2 {\let\theaccent#1\skewfactor#2
\mathpalette\putaccentunder}
\def\putaccentunder#1#2{\oalign{$#1#2$\crcr\hidewidth
\vbox to.2ex{\hbox{$#1\skew\skewfactor\theaccent{}$}\vss}\hidewidth}}
\def\name{\mathunderaccent\tilde-3 }

\newenvironment{PROOF}[2][\proofname.]
   {\begin{proof}[#1]\renewcommand{\qedsymbol}{\textsquare$_{\rm #2}$}}
   {\end{proof}}

\usepackage[hidelinks]{hyperref}

\begin{document}
\makeatletter\def\shfiuwefootnote{\gdef\@thefnmark{}\@footnotetext}\makeatother\shfiuwefootnote{Version 2022-08-12\_2. See \url{https://shelah.logic.at/papers/945/} for possible updates.}

\title{On CON($\gd_\lambda >$ cov$_\lambda$(meagre)) \\
Sh:945}

\author{Saharon Shelah}

\address{Einstein Institute of c\\
Edmond J. Safra Campus, Givat Ram\\
The Hebrew University of Jerusalem\\
Jerusalem, 91904, Israel\\
and \\
Department of Mathematics\\
Hill Center - Busch Campus \\ 
Rutgers, The State University of New Jersey \\
110 Frelinghuysen Road \\
Piscataway, NJ 08854-8019 USA}

\email{shelah@math.huji.ac.il}

\urladdr{http://shelah.logic.at}

\thanks{In earlier versions (up to 2019), 
the author thanks Alice Leonhardt for the beautiful typing.  This research was supported by the United States-Israel Binational Science Foundation.  First version (in TeX, September 16, 2005 as F741). First version in LaTeX January 21, 2009.}

\thanks{In later versions, the author would like to thank the typist for his work and is also grateful for the generous funding of typing services donated by a person who wishes to remain anonymous.}




\date{August 11, 2022}

\begin{abstract}
    We prove the consistency of: for suitable strongly
    inaccessible cardinal $\lambda$ the dominating number, i.e., the cofinality of ${}^\lambda \lambda$, is strictly bigger than cov$_\lambda$(meagre), i.e. the minimal number of nowhere dense subsets of ${}^\lambda 2$ needed to cover it.  This answers a question of Matet.  
\end{abstract}

\maketitle

\numberwithin{equation}{section}

\setcounter{section}{-1}

\newpage

\section{Introduction}\label{0}

Cardinal characteristics were defined, historically, over the continuum. See the celebrated Van Dowen \cite{VD}, for the general topologist perspective and  the 
excellent survey Blass \cite{Bls10}, Bartoszy\'nski \cite{Bar10} for the set theoretic perspective.  In
recent years there are many results concerning generalized cardinal characteristics.  The idea is to imitate the definition of a given characteristic over the continuum, by translating it to uncountable
cardinals.

It is reasonable to distinguish regular cardinals and singular cardinals.  Among the regular cardinals, it makes  sense to distinguish limit cardinals from successor cardinals.  In this paper we focus on
strongly inaccessible cardinals.  These cardinals and their characteristics behave, in many cases, much like $\aleph_0$, but certainly not always.  See Landver \cite{Land92}, Cummings-Shelah \cite{Sh:541} and Matet-Shelah \cite{Sh:804}. Our main result is
the consistency of $\cov_\lambda(\meagre) < \gd_\lambda$ at a super-compact cardinal $\lambda$, and we begin with the following definitions:

We shall define three cardinal invariants (but the paper deals, actually, just with two of them): 

\begin{definition}\label{z1}  
    The \emph{bounding} and \emph{dominating numbers}.  
    
    Let $\lambda$ be an inaccessible cardinal. For $f,g \in {}^\lambda \lambda$  let: 
    
    \begin{enumerate} 
        \item[(a)]  \emph{$f \le^* g$} if $|\{\alpha < \lambda:f(\alpha) > g(\alpha)\}| < \lambda$, 
        
        \item[(b)] $A \subseteq {}^\lambda \lambda$ is \emph{unbounded} if there is no $h \in {}^\lambda \lambda$ such that $f \in A \Rightarrow f \le^* h$, 
        
        \item[(c)]  $A \subseteq {}^\lambda \lambda$ is \emph{dominating} when for every $f \in {}^\lambda \lambda$ there exists $g \in A$ such that $f
        \le^* g$, 
        
        \item[(d)]  \emph{the bounding number} for $\lambda$, denoted by $\gb_\lambda$, is min$\{|A|:A \subseteq {}^{\lambda} \lambda$ is unbounded in ${}^\lambda \lambda\}$, 
        
        \item[(e)]  \emph{the dominating number} for $\lambda$, denoted by $\gd_\lambda$, is $ \min \{|A|:A \subseteq {}^{\lambda} \lambda$ is dominating in ${}^\lambda
        \lambda\}$.
    \end{enumerate}
    
    Notice that the usual definitions of $\gb$ and $\gd$ are $\gb_{\aleph_0}$ and $\gd_{\aleph_0}$ according to Definition \ref{z1}.  The definition of $\cov_\lambda$(meagre) involves some topology.  
\end{definition}

\begin{definition}\label{z5}
    \emph{The meagre covering number}. 
    
    Let $\lambda$ be a regular cardinal. 
    
    \begin{enumerate}
        \item[(a)]  ${}^\lambda 2$ is the space of functions from $\lambda$ into 2,  
        
        \item[(b)]  $({}^\lambda 2)^{[\nu]} := \{\eta \in {}^\lambda 2:\nu \triangleleft \eta\}$ and for $\nu \in {}^{\lambda >}2 := \bigcup\limits_{\alpha < \lambda} {}^\alpha 2$, 
        
        \item[(c)]  $\cU \subseteq {}^\lambda 2$ is \emph{open} in the topology $({}^\lambda 2)_{< \lambda}$, iff for every $\eta \in \cU,$ there exists $i < \lambda$ such that $({}^\lambda 2)^{[\eta \rest i]} \subseteq \cU$, 
        
        \item[(d)] $\cU \subseteq {}^{\lambda} 2$ is \emph{meagre} iff it is the union of $\leq \lambda$ no-where dense subsets,
        
        \item[(e)] $\cov_\lambda$(meagre) is the minimal cardinality of a  family of meagre subsets of $({}^\lambda 2)_{< \lambda}$ which covers this space.
    \end{enumerate}
\end{definition}

This paper deals with the relationship between $\gd_\lambda$ and $\cov_\lambda$(meagre).  If $\lambda$ is a successor cardinal then $\cov_\lambda(\meagre) < \gd_\lambda$ is consistent (see (b) below). Matet asked (a personal communication) whether $\gd_\lambda \le \text{ cov}_\lambda$(meagre) is provable in ZFC, where $\lambda$ is strongly inaccessible. We give here a negative answer.
    
For $\lambda$ a super-compact cardinal and $\lambda < \kappa = \text{cf}(\kappa) < \mu = \mu^\lambda$, we force large $\gd_\lambda$ i.e., $\gd_\lambda = \mu$ and small covering number (i.e., $\cov_\lambda$(meagre) $= \kappa$).  A similar result should hold also for a wider class of cardinals and we intend to return elsewhere to this subject.

Let us sketch some known results.  These results are related to the unequality number and the covering number for category.  Recall:

\begin{definition}\label{z14}
    \emph{The unequality number}.
    
    Let $\kappa$ be an infinite cardinal. The \emph{unequality number of $\kappa,{\mathfrak e}_\kappa$}, is the minimal cardinal $\lambda$ satisfying that there is a set $\cF \subseteq {}^\kappa \kappa$ of cardinality $\lambda,$ such that there is no $g \in {}^\kappa \kappa$ satisfying $(\forall f \in \cF)(\exists^\kappa \alpha < \kappa)(f(\alpha) = g(\alpha))$. 
    
    For $\kappa = \aleph_0,{\mathfrak e}_\kappa = 
    \cov_{\aleph_0}(\text{meagre})$; see Bartoszy\'nski (in \cite{Bar87}) and Miller (in \cite{Mi82}).  
    
    Now,
    
    \begin{enumerate}
        \item[(a)]  the statement ${\mathfrak e}_\kappa = \cov_\kappa(\text{meagre})$ is valid for $\kappa > \aleph_0$, in the case that $\kappa$ is strongly inaccessible, by \cite{Land92}.  But if $\kappa$ is a successor cardinal, it may fail,   
        
        \item[(b)]  if $\kappa < \kappa^{< \kappa}$ then $\cov_\kappa(\meagre) = \kappa^+$.  This is due to Landver (in \cite{Land92}).
    \end{enumerate}
\end{definition}

We intend also to address:

\begin{problem}\label{z15}
    Can we replace ``super-compact" by ``strongly inaccessible"?
\end{problem}

\begin{problem}\label{z16}\ 

    (1) Can we prove the consistency of $\cov_\lambda(\meagre) < \gb_\lambda$?

    (2) For $\lambda$ strongly inaccessible (or just Laver indestructible super-compact) is there a non-trivial $\lambda^+$-c.c. $(< \lambda)$-strategically complete forcing notion $\bbQ$ which is ${}^\lambda \lambda$-bounding?
       
    We say more in subsequent works \cite{Sh:1004}, \cite{Sh:E82} and in preparation \cite{Sh:1100}.
    
    A point which in a previous version was just a step along the way, the referee asked to justify fully, was analyzed to be serious. This was done but eventually is separated to \cite{Sh:1126}.  
    A posteriori the point is that in the parallel case for  $\lambda = \aleph_0$, for full memory FS iteration  such a claim is true.  In fact, 
    by Judah-Shelah \cite{Sh:292}, if $\langle \bbP_\alpha, \name{\bbQ}_\beta:\alpha \le \alpha(*),\beta <  \alpha(*)\rangle$ is FS iteration of Suslin-c.c.c. forcing notion, $\name{\bbQ}_\beta$ with the generic $\name\eta_\beta \in {}^\omega \omega$ and for notational transparency, its definition is with no parameter and $\zeta:\beta(*) \rightarrow \alpha(*)$ is increasing
    and $\bbP = \langle \bbP'_\alpha,\name{\bbQ}'_\beta:\alpha \le
    \beta(*),\beta < \beta(*)\rangle$ is FS iteration, but $\name{\bbQ}'_\beta$ is defined exactly as
    $\name{\bbQ}_{\zeta(\beta)}$ is but now in $\mathbf V^{\bbP'_\beta}$ rather then in $\mathbf V^{\bbP_{\zeta(\beta)}}$ then $\Vdash_{\bbP_{\alpha(*)}} ``\langle \name\eta_{\zeta(\beta)}:\beta < \beta(*)\rangle$ is
    generic for $\bbP'_{\beta(*)}$ over $\mathbf V"$.
    
    Now this is not clear to us for $(< \lambda)$-support iteration of $(< \lambda)$-strategically complete forcing notions  for $ \lambda > {\aleph_0} $.   
    The solution is essentially to change the iteration: to use a ``quite generic" $(< \lambda)$-support iteration which ``includes" the one we like and use
    the complete sub-forcing it generates; see \cite{Sh:1126}.  
\end{problem}

We thank the referee, Shimoni Garti and Haim Horowitz for helpful comments and pressuring me to expand some proofs and  Johannes Sch\"urz and Martin Goldstern for pointing  several times problem in \cite{Sh:945}, in particular point out in 2019 that in an earlier version  of the proof of \cite[2.7=La32]{Sh:945} the statement $\circledast'_4$ was insufficient and  later pointing out a problem in earlier version of the end of the proof of \cite[3.43=Le70]{Sh:1126} which motivate the addition of $\Vdash_{\bbP_{\bfm}}$``$\{ \name{\eta}_{s}: s \in M_{\bfn}\}$ is cofinal in $\left( \Pi_{\varp < \lambda \theta_{\varp}}, <_{J_{\lambda}^{\rm{bd}}} \right)$'' and being reasonable. 
    
We try to use standard notation.  We use $\theta,\kappa,\lambda,\mu,\chi$ for cardinals and $\alpha,\beta,\gamma,\delta,\varepsilon,\zeta$ for ordinals.  We use also $i$ and $j$ as ordinals.  We adopt the Cohen convention that $p \le q$ means that $q$ gives more information, in forcing notions.  The symbol $\triangleleft$ is preserved for ``being an initial segment". Also recall ${}^B A = \{f:f$ a function from $B$ to $A\}$ and let ${}^{\alpha >}A = \cup\{{}^\beta A:\beta < \alpha\}$, some prefer ${}^{<\alpha} A$, but ${}^{\alpha >} A$ is used systematically in the author's papers.  Lastly, $J^{\text{\rm bd}}_\lambda$ denotes the ideal of the bounded subsets of $\lambda$. 
    
For exact references to \cite{Sh:1126} see the introduction there, just before \cite[0.1= Lz6]{Sh:1126}. 
The picture of cardinal invariants related to uncountable $\lambda$ is related but usually quite different than the one for $\aleph_0$, they are more similar if $\kappa$ is ``large" enough, mainly strongly inaccessible.

\newpage

\section{Preliminaries}\label{1}

\begin{definition}\label{z17}
    Let $\lambda$ be super-compact.  We say that \emph{$h:\lambda \rightarrow  \cH(\lambda)$ is a Laver diamond (for $\lambda$)} \when \, for every 
    $x \in \mathbf V$ there are a normal fine ultrafilter $D$  over $I = [\cH(\chi)]^{<\lambda}$ for some $\chi,$ such that $x \in \cH(\chi)$ and the Mostowski collapse $\mathbf j$ of $\mathbf V^I/D$
    maps $\langle h(\sup(u \cap \lambda)):u \in I\rangle/D$ to $x$ (we can use elementary embeddings instead of an ultrafilter).
\end{definition}

\begin{notation}\label{z19}
    If $\bbP$ is a forcing notion in $\mathbf V$ then $\mathbf V^{\bbP}$ denotes $\mathbf V[\mathbf G]$ for $\mathbf G \subseteq \bbP$ generic over $\mathbf V$; we may write $\mathbf V[\bbP]$ instead.
\end{notation}

The most straightforward way to increase $\gb_\lambda$ in the classical case of $\aleph_0$ is Hechler forcing (dominating real forcing).  A condition is a function $f_p:\omega \rightarrow \omega$ which is separated into a finite trunk $\eta_p$ and the rest of the function. Formally, $p = (\eta_p,f_p)$ where $\eta_p \trianglelefteq f_p$.
    
If $p,q$ are conditions then $p \le q$ \Iff \, $\eta_p \trianglelefteq \eta_q$ and $f_q(n) \ge f_p(n)$ for every $n \notin \dom(\eta_p)$ hence for every $n$.  A generic object adds a function  $g:\omega \rightarrow \omega$ which dominates the functions from the ground model.  By iterating Hechler reals one increases the bounding number $\gb$.
    
If $\lambda = \lambda^{< \lambda}$ then one can define the \emph{generalized Hechler forcing} $\bbD_\lambda$ by replacing $\omega$ by $\lambda$. The basic step is $(< \lambda)$-complete and $\lambda^+$-c.c. and actually
$\lambda$-centered. Hence one can iterate and increase $\gb_\lambda$.
    
In \cite[\S1,\S2]{Sh:326} and then Goldstern-Shelah \cite{Sh:448}, Kellner-Shelah \cite{Sh:961}  other invariants are discussed.  Consider two functions $f,g:\omega \rightarrow (\omega \backslash \{0\})$ going to infinity such that $f \ge g$ and ask about:
    
\begin{itemize}
    \item  $\gc^+_{f,g} = \min\{|\cF|:\cF \subseteq \prod\limits_{i} [f(i)]^{g(i)}$ and $(\forall \eta \in \prod\limits_{i} f(i))(\exists g \in \cF)[\bigwedge\limits_{i} \eta(i) \in g(i)] \}$, 
    
    \item  $\gc^-_{f,g} = \min\{\cF:\cF \subseteq \prod\limits_{i} f(i)$ and for no $g \in \prod\limits_{i} [f(i)]^{g(i)}$ do we have $(\forall \eta \in \cF)(\forall^\infty i)(\eta(i) \in g(i)) \}.$
\end{itemize}
    
There are relevant forcing notions; we could have used $[f(i)]^{<g(i)},$ this case is generalized here, so below $f = gu$ replaced by $\langle \theta_{\varp}: \varp < \lambda \rangle$ we shall use a $\lambda^+$-c.c. one as in c.c.c. creature forcing (see more in \cite{Sh:628},\cite{Sh:1067}). 

For transparency,

\begin{convention} \label{z21}
    Below $\lambda,\bar\theta$ are as in \ref{z23} below.
\end{convention}

\begin{definition}\label{z23}
    Let $\lambda$ be inaccessible, $\bar\theta = 
    \langle \theta_\varepsilon:\varepsilon < \lambda\rangle$ be a sequence of regular cardinals $< \lambda$ satisfying $\theta_\varepsilon > \varepsilon$.
    
    (1) We define the forcing notion $\bbQ = \bbQ_{\bar\theta}$ by:  
    
    \begin{enumerate}
        \item[(A)]  $p \in \bbQ$ \underline{iff}: 
    
        \begin{enumerate}
            \item[(a)] $p = (\eta,f) = (\eta^p,f^p)$,   
        
            \item[(b)]  $\eta \in \prod\limits_{\zeta < \varepsilon} \theta_\zeta$ for some $\varepsilon < \lambda$, ($\eta$ is called the trunk of $p$), 
        
            \item[(c)] $f \in \prod\limits_{\zeta < \lambda} \theta_\zeta$, 
        
            \item[(d)] $\eta \triangleleft f$.
        \end{enumerate}
        
        \item[(B)]  $p \le_{\mathbb Q} q$ iff:  
        
        \begin{enumerate}
            \item[(a)]  $\eta^p \trianglelefteq \eta^q$,  
            \item[(b)]  $f^p \le f^q$, i.e. $(\forall \varepsilon < \lambda)
            f^p(\varepsilon) \le f^q(\varepsilon)$, 
            
            \item[(c)]  if $\ell g(\eta^p) \le \varepsilon < \ell g(\eta^q)$ then $\eta^q(\varepsilon) \in [f^p(\varepsilon),\lambda)$, actually follows.
        \end{enumerate}
    \end{enumerate}
    
    (2) The generic is $\name\eta = \bigcup \, \{\eta^p:
    p \in \name{\mathbf G}_{\bbQ_{\bar\theta}}\}$.
\end{definition}

The new forcing defined above is not $(< \lambda)$-complete anymore.  By  fixing a trunk $\eta \in \Pi_{\varp < \zeta} \theta_{\varp}$ one can define a short increasing sequence $\bar{p} = \langle p_{\varp}: \varp < \zeta \rangle$ of conditions which goes up to $\theta_\zeta$ at the $\zeta$-th coordinate that is, $\langle p_{\varp}(\zeta): \zeta < \theta_{\zeta} \rangle$ is increasing and hence has no upper bound in $\theta_\zeta,$ so $\bar{p}$ has no upper bound in $\bbQ$.  However, this forcing is $(< \lambda)$-strategically complete (see Definition \ref{z28}(2), below) since the COM (= completeness) player can increase the trunk at each move.

\begin{remark}\label{z25}
    The forcing is parallel to the creature forcing from
    \cite[\S1,\S2]{Sh:326}, \cite{Sh:961} but they are ${}^\omega \omega$-bounding.
\end{remark}

Recall,

\begin{definition}\label{z28}\ 

    (1) We say that a forcing notion
    $\mathbb P$ is \emph{$\alpha$-strategically complete} \underline{when}  for each $p \in \mathbb P$ in the 
    following game $\Game_\alpha(p,\mathbb P)$ between the players COM and INC, the player COM has a winning strategy.
    
    A play lasts $\alpha$ moves; in the $\beta$-th move, first the player COM chooses $p_\beta \in \mathbb P$ such that $p \le_{\mathbb P} p_\beta$ and $\gamma < \beta \Rightarrow q_\gamma \le_{\mathbb P} p_\beta$ and second the player INC chooses $q_\beta \in \mathbb P$ such that $p_\beta \le_{\mathbb P} q_\beta$.
    
    The player COM \emph{wins} a play if he has a legal move for every $\beta < \alpha$.
    
    (2) We say that a forcing notion $\mathbb P$ is \emph{$(< \lambda)$-strategically complete} \underline{when} it is $\alpha$-strategically complete for every $\alpha < \lambda$.
\end{definition}

Basic properties of $\bbQ_{\bar\theta}$ are summarized and proved in \cite[\S2]{Sh:949}. 

The following fact describes some immediate connections between various concepts of completeness:

\begin{fact}\label{z30}\ 

    \begin{enumerate}
        \item[(a)]  If $\bar{\bbQ} = 
        \langle \bbP_\alpha,\name{\bbQ}_\beta:\alpha \le \delta,\beta < \delta\rangle$ is a $(< \lambda)$-support iteration of $(< \lambda)$-strategically complete forcing notions, \then \, $\bbP_\delta$ is also $(< \lambda)$-strategically complete (see e.g. \cite{Sh:546}),
        
        \item[(b)]  If $\bbP$ is $(< \lambda)$-strategically complete forcing notion \underline{then} $({}^{\lambda >} \Ord)^{\mathbf V} = ({}^{\lambda >} \Ord)^{\mathbf V^{\bbP}}$, and consequently $\lambda$ is strongly inaccessible in $\mathbf V^{\bbP}$,
        
        \item[(c)] like (a) replacing $\lqq (< \lambda)$-strategically complete" by ``$(< \lambda)$-complete"or by ``$\alpha$-strategically complete'',
          
        \item[(d)]  if $\bbP$ is $(< \lambda)$-complete then $\bbP$ is $(< \lambda)$-strategically complete.  
    \end{enumerate}
\end{fact}

\begin{definition}\label{z32}
    For an ordinal $\alpha_* := \alpha(*)$ let $\mathbf
    Q_{\lambda,\bar\theta,\alpha(*)}$ be the class of quintuple  $\mathbf q = (\bar u, \bar{ {\mathscr P}}, 
    \bar{\bbP},\name{\bar{\bbQ}},\name{\bar\eta})$
    consisting of (omitting $\alpha_*$ means for some $\alpha_*$ and  we let $\ell g(\mathbf q) = \alpha_{\mathbf q} = \alpha_*$):
    
    \begin{enumerate}
        \item[(a)] $\bar u = \langle u_\alpha:\alpha < \alpha_*\rangle$ and $\bar{\cP} = \langle \cP_\alpha:\alpha < \alpha_*\rangle$ where $\cP_\alpha \subseteq [u_\alpha]^{\le\lambda},u_\alpha \subseteq \alpha$, \wilog \, $\cP_\alpha$ is closed under subsets (but is not necessarily an ideal), 
        
        \item[(b)]  $\langle \bbP_{0,\alpha},\name{\bbQ}_{0,\beta}: \alpha \le \alpha_*,\beta < \alpha_*\rangle$ is a  $(< \lambda)$-support iteration and let $\bbP_{\mathbf q,0} = \bbP_{\mathbf q,0,\alpha(\mathbf q)};$ but we may write $\bbP_{\alpha}, \name{\bbQ}_{\beta},$ 
        
        \item[(c)]  each $\bbP_\alpha$ is   $(<\lambda)$-strategically complete and $\lambda^+$-c.c.,
        
        \item[(d)]  $\name\eta_\beta \in \prod\limits_{\varepsilon < \lambda} \theta _ \varepsilon $  is the generic of $\name{\bbQ}_\beta$ where $\name\eta_\beta$, the generic of $\name{\bbQ}_\beta$ (defined in clause (e) below) is $\cup\{\eta_p:p \in \mathbf G_{\bbQ_\beta}\}$,  
        
        \item[(e)]  if $ \alpha < \alpha _* $  
        $\mathbf G \subseteq \bbP_\beta$ is generic over $\mathbf V$ then $\name\eta_\alpha[\mathbf G]$ in  $( \prod\limits_{\varepsilon < \lambda} \theta _  \varepsilon  ,<_{J^{\bd}_\lambda})$ dominates every $\nu \in  \prod\limits_{\varepsilon < \lambda} \theta _ \varepsilon$ from $\mathbf V[\langle \name\eta_\gamma:\gamma \in  u\rangle]$ when $u \in \cP_\alpha$; moreover, in $\mathbf V[\mathbf G]$:  
        
        \begin{enumerate}
            \item[$(*)$]  $\name{\bbQ}_\beta[\mathbf G]$ is the sub-forcing of $\bbQ_{\bar\theta}$  consisting of the $p \in \bbQ_{\bar\theta}$ such that: for some $\bar s,\name f,\eta_p$ (so $\eta_p = \eta$, etc.) we have:
            
            \begin{enumerate}
                \item[$(\alpha)$]  $p = (\eta,f) = (\eta_p,f_p)$ so $\eta \in \prod\limits_{\varp < \zeta} \theta_\varp$ for some $\zeta < \lambda$,  
                
                \item[$(\beta)$]  $\bar s = \langle (u_i,f_i):i < i_*\rangle$,  
                
                \item[$(\gamma)$]  $i_* < \lambda$, 
                
                \item[$(\delta)$]  for each $i < i_*$ we have $u_i \in \cP_\beta, \eta \triangleleft f_i \in \prod\limits_{\varepsilon < \lambda} \theta_ \varepsilon 
                $ and $f_i \in \mathbf V[\langle \name\eta_\gamma[\mathbf G]: \gamma \in u_i\rangle]$, 
                
                \item[$(\varp)$]  $f = \sup\{f_i:i < i_*\}$, i.e. $\varp < \lambda \Rightarrow f(\varp) = \cup\{f_i(\varp):i < i_*\}$;we may add $ i_* < \theta _{\lg(\eta )} $.  
            \end{enumerate}
        \end{enumerate}
        
        \item[(f)] \underline{notation}: so $u_{\bfq,\alpha} = u_\alpha,\bbP_{\bfq,\alpha} = \bbP_\alpha$, etc., but when $\bfq$ is clear from the context we may omit it.
    \end{enumerate}
\end{definition}

\begin{definition}\label{z33}
    For $\bfq \in \bfQ_{\lambda,\bar\theta,\alpha(*)},$
    
    (1) We let $\alpha \le \alpha_*,\bbP_{1,\alpha} = \bbP^{\bfq}_{1,\alpha}$ be essentially the completion of $\bbP_\alpha$; we express it by:
    
    \begin{enumerate}
        \item[$(*)_1$]  the elements of $\bbP_{1,\alpha} = \bbP_{\bfq,1,\alpha}$ are of the form  $\mathbf B(\ldots,\name\eta_{\gamma_i},\ldots)_{i < i(*)},$ where:
        
        \begin{enumerate}
            \item[$(\alpha)$]  $i_* = i(*) \le \lambda$, 
            
            \item[$(\beta)$]  $\gamma_i  < \alpha$ for $i < i_*$,  
            
            \item[$(\gamma)$]  $\mathbf B$ is a $\lambda$-Borel function from ${}^{i(*)}(\prod\limits_{\varepsilon < \lambda}  \theta_ \varepsilon)$ into $\{0,1\} = \{\false,\true\}$; $\mathbf
            B$ is from $\mathbf V$, of course, such that $\nVdash_{\bbP_{\mathbf q}} ``\mathbf B(\ldots,\name\eta_{\gamma_i},\ldots)_{i<i(*)} = 0"$. 
        \end{enumerate}
        
        \item[$(*)_2$]  the order is natural: $\bbP_{1,\alpha} \models ``\mathbf B_1(\ldots,\name\eta_{\gamma(i,1)},\ldots)_{i<i(1)} \le  \mathbf B_2(\ldots,\name\eta_{\gamma(i,2)},\ldots)_{i<i(2)}"$ iff $\Vdash_{\bbP_\alpha}$ ``if 
        $\mathbf B_2(\ldots,\name\eta_{\gamma(i,2)}[\name{\mathbf G}],\ldots)_{i<i(2)}$ is equal to 1 then so is 
        $\mathbf B_1(\ldots,\name\eta_{\gamma(i,1)},\ldots)_{i<i(1)}"$.
    \end{enumerate}
    
    (2) For $\cU \subseteq \alpha_*$ let $\bbP_{\cU} = \bbP^{\bfq}_{\cU}$ be the sub-forcing of $\bbP_{1,\alpha(\mathbf q)}$ consists  of  $\{\mathbf B(\ldots,\name\eta_{\gamma(i)},\ldots)_{i<i(*)} \in  \bbP_{1, \alpha(\mathbf q)}:i(*) \le \lambda$ and $\gamma_i \in \cU$ for every $i <i(*)\}$.
\end{definition}

\begin{claim}\label{z35}\ 

    (1) For any sequence $\langle u_\alpha,\cP_\alpha:\alpha < \alpha_*\rangle$ as above, i.e. as in clause (a) of Definition \ref{z32}, there is one and only one $\mathbf q \in \bfQ_{\lambda,\bar\theta,\alpha_*}$ with $u_{\bfq,\alpha} = u_\alpha,\cP_{\bfq,\alpha} = \cP_\alpha$ for $\alpha < \alpha_*$.
    
    (1A) For $\alpha \le \alpha_*$,  the forcing notions
    $\bbP_{\bfq,1,\alpha},\bbP_{\mathbf q,1,\cU}$'s are well defined and are as demanded in Definition \ref{z33}.
    
    (2) For every  $ \alpha \le \alpha_*$ $\bbP^\bullet_{\bfq,\alpha} = \bbP_{\bfq, 0,\alpha}^{\bullet}$ is $\bbP_{\bfq, \alpha} = \bbP_{\bfq, 0, \alpha}$ restricted to the set of $p \in \bbP_{\bfq, \alpha}$  (from Definition \ref{z32}) satisfying the following is dense in $\bbP_{\bfq, \alpha}$: 
    
    \begin{enumerate}
        \item[$(*)$]  if $\beta \in \dom(p)$, then $q = p(\beta)$ is a $\bbP_\beta$-name of a member of $\name{\bbQ}_\beta$ such that:
        
        \begin{enumerate}
            \item[(a)]  $\eta_q,i_q, \langle u_{q,i}:i < i_q  \rangle$ are objects (not just $\bbP_\beta $-names),  
            
            \item[(b)]  $\name f_q = \sup\{\name f_i:i < i_q\}$, each $\name f_i$ is a $\bbP_\beta$-name of a member of $\prod\limits_{\varepsilon < \lambda} \theta _ \varepsilon $,  
            
            \item[(c)]  each $\name f_i$ has the form $\mathbf B_{q,i}
            (\ldots,\name\eta_{\gamma(i,j)},\ldots)_{j<j(*) \le \lambda}$ where $\{\gamma(i,j):j<j(*)\} \subseteq u_{q, i } $  and $\bfB_q$ is a Borel function from ${}^{   j(*)}(\prod\limits_{\varepsilon < \lambda}  \theta _  \varepsilon  )$ into $\prod\limits_{\varepsilon < \lambda} \theta_ \varepsilon,$  
            
            \item[(d)] $p(\beta) = (\name\eta_q,\name f_q)$.
        \end{enumerate}
    \end{enumerate}
    
    (2A) Abusing our notation we may identify $\bbP_{\alpha}^{\bullet}$ with $\bbP_{\alpha}.$ 
    
    (3) Above for every $v \subseteq \alpha$ and $j_* < \lambda$ the set of $p \in \bbP^\bullet_\alpha$ such that $v \subseteq \dom(p) \wedge (\forall \beta \in \dom(p))(\ell g(\eta_{p(\beta)}) > j_*)$ is dense.
    
    (4) $\bbP^{\bullet}_{\mathbf q,0,\alpha} \lessdot \bbP_{\mathbf q,1,\alpha}$ moreover $\bbP^{\bullet}_{\mathbf q,0,\alpha}$ is dense in $\bbP_{\mathbf q,1,\alpha}$ and in $\bbP_{\bfq, 0, \alpha}$ and $\cU_1 \subseteq \cU_2 \subseteq \alpha \leq \alpha_{\mathbf q} \Rightarrow \bbP_{\bfq, 1, \cU_{q}} \lessdot \bbP_{\bfq, 1, \cU_{2}}
    \lessdot \bbP_{\mathbf q,1,\alpha}$ so $\bbP_{\mathbf q,1,\{\beta:\beta < \alpha\}} = \bbP_{\mathbf q, 1,
    \alpha}$ and $|\bbP_{\mathbf q,1,\cU}| \le |\cU|^\lambda$.
    
    (5) If $\alpha < \alpha_*$ and $u \in \cP_\alpha$ then
    $\name\eta_\alpha \in \prod\limits_{\varepsilon < \lambda} \theta _ \varepsilon  $ dominates every $\nu \in (
    \prod\limits_{\varepsilon < \lambda} \theta
    _\varepsilon)^{\mathbf V[\name{\bar\eta} \rest u]}$ and $\bfV[\bfG_{\bbP_{\bfq, 1, \cU}}] = \bfV[\langle \eta_{\alpha}: \alpha \in \cU \rangle],$ where $\eta_{\alpha} = \name{\eta}[\bfG_{\bbP_{\bfq}, 1, \alpha}].$
    
    (6) Assume $\mathbf G \subseteq \bbP_{\mathbf q}$ is generic over $\mathbf V,\name\eta_\alpha = \name\eta_\alpha[\mathbf G]$ and $\eta'_\alpha \in
    (\prod\limits_{\varepsilon < \lambda } 
    \theta_ \varepsilon )^ {\mathbf V[\mathbf G]}$ for $\alpha < \alpha_*$ and $\{(\alpha,\varp):\alpha< \alpha_*,\varp < \alpha$ and $\eta_\alpha(\varp) \ne \eta'_\alpha(\varp)\}$ has cardinality $< \lambda$.  \Then \, for some (really unique) $\mathbf G'$ we have $\mathbf G' \subseteq \bbP_{\mathbf q}$ is generic over $\mathbf V$ and $ \mathbf{V} [ \mathbf{G} ']= \mathbf{V} [\mathbf{G} ]$ and $ \name{ \eta } _ \alpha [\mathbf{G} ]= \eta ' _ \alpha $ for $ \alpha < \alpha _*.$ 
    
    (7) Like (6) for $ {\bbP}^{\mathbf q}_{{\cU}}.$
\end{claim}

\begin{PROOF}{\ref{z35}}
    See \cite[1.13=Lc8, 1.16=Lc11]{Sh:1126}.
\end{PROOF}

\begin{theorem}\label{z38}
    For any ordinal $\alpha_*$ there is a quadruple 
    $(\mathbf q,\delta_*,\cU_*,h)$ such that:
    
    \begin{enumerate}
        \item[(A)]
        
        \begin{enumerate}
            \item[(a)]  $\mathbf q \in \mathbf Q_{\lambda,\bar\theta}$ and let $\delta_* = \ell g(\mathbf q),$
            
            \item[(b)]  $\cU_* \subseteq \delta_*$ has order type $\alpha_*,$
            
            \item[(c)]  $h$ is the order preserving function from $\alpha_*$ onto $\cU_*,$
            
            \item[(d)]  if $\alpha \in \cU_*$ then  $\cU_* \cap \alpha \in \cP_{\bfq,\alpha},$ 
            
            \item[(e)]      if $ \cf(\alpha _*) >  \lambda $ then  in $ \mathbf{V} ^{\mathbb{P} _ \mathbf{q} }$ the set $ \{\name{ \eta } _ \alpha : \alpha \in {\mathscr U } _ * \} $ is cofinal in $ (\Pi _{\varepsilon < \lambda }\theta _ \varepsilon, \le _{J^{\bd}_\lambda}).$  
        \end{enumerate}
        
        \item[(B)]  if $\cU_1 \subseteq \cU_*,\cU_2 \subseteq \cU_*, \otp(\cU_1) = \otp(\cU_2)$ and $g$ is the order preserving function from $\cU_1$ onto $\cU_2$, \then \, $g$ induces an isomorphism $\hat g$ from $\bbP_{\mathbf q,\cU_1}$ onto $\bbP_{\mathbf q,\cU_2}$ mapping $\name\eta_\beta$ to $\name\eta_{g(\beta)}$ for $\beta \in \cU_1$.
    \end{enumerate}
\end{theorem}

\begin{PROOF}{\ref{z38}}
    By \cite[2.14=Le38]{Sh:1126}, in particular  clause (A)(e) is justified by clause (E) there.
\end{PROOF}

\begin{definition}\label{z41}\

    (1) Let ``$\bbP$ is a $(< \lambda)$-\dic" mean: 
    
    \begin{enumerate}
        \item[$(*)$]  if $J$ is a directed partial order of cardinality $< \lambda$ and $p_s \in \bbP$ for $s \in J$ and $s \le_J t \Rightarrow p_s \le_{\bbP} p_t$ \then \, the set $\{p_s:s \in J\}$ has an upper bound in $\bbP$.
    \end{enumerate}
    
    (2) Assume $ \bar{ \theta  } = \langle \theta _\varepsilon: \varepsilon < \lambda \rangle$ and $ \theta _ \varepsilon = \cf(\theta_\varepsilon ) \in (\varepsilon, \lambda )$.  We say  ``$\bbP$ is $\dc $-directed complete" \when:
    
    \begin{enumerate} 
        \item[($*$)]  if $\mathbb{P} \in N \prec ({\mathscr H} (\chi ), \in )$ and $ \lambda _N = N \cap \lambda $ is inaccessible $ < \lambda$ and $ N^{< \lambda _N}  \subseteq N$, $ \| N \cap \theta _* \| < \theta _{\lambda_N}$ and $\mathbf{G} \subseteq  \mathbb{P} \cap N $ is $ (N, \mathbb{P} )$-generic \then \, $ \mathbf{G} $ has a common upper bound
    \end{enumerate} 
\end{definition} 

\begin{definition}\label{z44}
    We say that $ h $ is a $\dc$-Laver diamond when:
    
        \begin{enumerate} 
            \item[(a)] $\lambda $ is a super-compact cardinal, 
            
            \item[(b)] $\theta _* > \lambda $, 
            
            \item[(c)] $\bar{ \theta } = \langle \theta _ \varepsilon: \varepsilon < \lambda \rangle,$ where $ \varepsilon < \theta _ \varepsilon = \cf( \theta _ \varepsilon ) < \lambda $,
            
            \item[(d)] for every  $ x $ and $ \chi > \lambda $ such that  $ x \in {\mathscr H} (\chi )$ there are $ \mathbf{j}, \mathbf{M}$ such that:
            
            \begin{enumerate} 
                \item[$ \bullet_1$] $\mathbf{M} $  is a transitive class, 
                
                \item[$ \bullet _2$] $ \mathbf{j} $ is an elementary embedding of $ \mathbf{V}$ into $ \mathbf{M},$
                
                \item[$ \bullet _3$] $ \mathbf{M} ^ {\le \chi }  \subseteq \mathbf{M} $, 
                
                \item[$ \bullet _4$] $ \mathbf{j} $ is the identity on $  {\mathscr H} (\lambda)$,
                
                \item[$ \bullet _5$] $ \mathbf{j}(h)(\lambda ) = x$ and $\mathbf{j} (\bar{ \theta }  )(\lambda ) = \theta _* $,  
            \end{enumerate} 
    \end{enumerate} 
\end{definition} 

\begin{observation}\label{z47}\ 

    (1) There are enough cases of $ h $ as in \ref{z44}, for example,  for every  Laver diamond $ h $, 
    
    \begin{enumerate} 
        \item[(a)] if $ \theta = (\beth_{\lambda ^+})^{+}$ and $\langle  \theta _ \varepsilon : \varepsilon < \lambda \rangle$ defined by  $ \theta _ \varepsilon = \beth (\| \varepsilon \|^+)^{+}$  \then \, $ h $ is a Laver diamond for $ (\lambda, \theta_ *, \bar{ \theta}),$    
        
        \item[(b)] there is $\bar{ \theta } $ such that for every regular  cardinal $ \theta _* >  \lambda $, $ h $ is a Laver diamond for $ (\lambda, \theta_ *, \bar{ \theta } )$.
    \end{enumerate} 
    
    (2) Being  a $ (< \lambda )$-directed complete forcing notion, is preserved by $ (< \lambda )$-support iteration. Similarly for being ``$(< \lambda, \theta_{\ast}, \bar{\theta})$-directed complete''. 
    
    (3) If the forcing notion $ \mathbb{Q} $ is $ (< \lambda  )$-\dic and $ \zdc $ from Definition \ref{a8}(2) below \then \, $  \mathbb{Q} $ is $ \dc $-\dic. 
\end{observation}

\begin{PROOF}{\ref{z47}}
    Easy, e.g for clause (b) in part (1), let $\bar{  \theta }  = \langle  \theta _ \varepsilon : \varepsilon < \lambda \rangle$ where for $ \varepsilon < \lambda $ we let $ \theta _ \varepsilon$ be the first regular cardinal $ \theta $ such that  $ \ge \varepsilon, $ and $ h(\varepsilon) = (\theta, \chi),$ for some $\chi$.   
\end{PROOF}

\newpage

\section{The forcing}\label{2} 

In this section we prove the main result of the paper, which reads as follows:

\begin{theorem}\label{a2} 
    Assume,
    
    \begin{enumerate}
        \item[(a)]  $\lambda$ is super-compact
        
        \item[(b)]   $\lambda < \kappa = \cf(\kappa) 
        < \mu = \cf(\mu) = \mu^\lambda$.
    \end{enumerate}
    
    \Then \, for some $(< \lambda)$-strategically complete $\lambda^{+}$-cc forcing notion hence $\bbP$ not collapsing cardinals $\ge \lambda,\lambda$ is still super-compact in $\mathbf V^{\bbP}$ and
    $\cov_\lambda${\rm (meagre)} $=\kappa,{\mathfrak d}_\lambda =\mu$.
\end{theorem}

\begin{PROOF}{\ref{a2}}
    Let $ \theta  _* = \cf(\theta) > 2^ \mu $ and $ \bar{\theta }$ be such that there is a Laver diamond for $ \dc $, justified by \ref{z47}.  By Lemma \ref{a20}(1)  below we  can (by \ref{a20}(2) below) 
    force $\zdc$ while maintaining the super-compactness of $\lambda$.  By Lemma \ref{a32} we force $\gd_\lambda = \mu \wedge \cov_\lambda(\meagre) = \kappa$ using a
    forcing notion $\bbP$ which is $\lambda^+$-c.c. and $(<
    \lambda)$-strategically complete.  Notice that $\lambda$ is still super-compact in the generic extension, so we are done.
\end{PROOF}

\begin{definition}\label{a8}\ 

    (1) For $\lambda$ super-compact we define $\boxdot_\lambda$ by:
    
    \begin{enumerate}
        \item[$\boxdot_\lambda$]  for any regular cardinal $\chi > \lambda$ and forcing notion $\bbP \in \cH(\chi)$ which is $(<\lambda)$-strategically complete (see Definition \ref{z28}(2)) the  following set ${\cS} = {\cS}_{\bbP} = \cS_{\chi,\bbP}$ is a  stationary subset of $[{\cH}(\chi)]^{< \lambda}$: ${\cS} = \cS_{\bbP} = \cS_{\chi,\bbP}$ is the set of $N$'s such that for some  $\lambda_N,\chi_N,\mathbf j = \mathbf j_N,\bbA = \bbA_N,M = M_N,\mathbf G = \mathbf G_N$ we have (and we may say $(\lambda_N,\chi_N,\mathbf j_N,\bbA_N,M_N,\mathbf G_N)$ is a $  \boxdot _ \lambda $-witness  for $N \in \cS_{\chi,\bbP}$ or for $(N,\bbP,\chi)$):  
        
        \begin{enumerate}
            \item [(a)]  $N \prec ({\cH}(\chi)^{\mathbf V},\in)$ and $\bbP \in N$, 
            
            \item[(b)]   the Mostowski collapse of $N$ is $\bbA$ and let $\mathbf j_N:N \rightarrow \bbA$ be the unique isomorphism, 
        
            \item[(c)]  $  N \cap \lambda$  is a  strongly inaccessible cardinal called $  \lambda _N = \lambda (N),$
             
            such that ${}^{\lambda(N)>}N
            \subseteq N$,  
            
            \item[(d)]  $\bbA \subseteq M := ({\cH}(\chi_N), \in), M$ is transitive as well as $\bbA$,  
            
            \item[(e)]  $\mathbf G \subseteq \mathbf j_N(\bbP)$ is generic over $\bbA$ for the forcing notion $\mathbf j_N(\bbP)$, 
            
            \item[(f)]  $M = \bbA[\mathbf G]$.
        \end{enumerate}
    \end{enumerate}
    
    (2) Assume  $ \lambda, \theta _*, \bar{ \theta } $ 
    satisfies clauses (b), (c) of Def \ref{z44} and $ \lambda$ is super-compact.  We define $ \boxdot_{< \lambda, \theta _*, \bar{ \theta }}$ by:

    \begin{enumerate} 
        \item[$\boxdot_{< \lambda, \theta _*, \bar{ \theta } }$] As in part (1) adding $ \| N \cap \theta _* \| <  \theta _{\lambda _N}$ and even  $ 2^{ \| N \cap \theta _*\| } <  \theta _{\lambda _N} $.   
    \end{enumerate} 
\end{definition} 

Our first lemma is closed to Laver's indestructibility.
It consists of two parts.  In the first part we prove that one can force $\boxdot_\lambda$ at a super-compact cardinal $\lambda$ while preserving its super-compactness.  In the second part, we deal with a more informative version $ \zdc $.  See \ref{z47},  this can be done in an indestructible manner.  Namely, any further extension of the universe by a  $ \zdc $-\dic forcing notion will preserve the $\dc$-Laver diamond and the principle  $\boxdot_{\lambda, \theta_{\ast}, \bar{\theta}}.$  

\begin{lemma}\label{a20}\ 

    (1) If $\lambda$ is super-compact (in the universe $ \mathbf{V} = \mathbf{V} _0$) \underline{then} after
    some preliminary forcing  $ \mathbb{R}$ of cardinality $\lambda$, getting a universe $ \mathbf{V}_1 = \mathbf{V} ^{\mathbb{R} }_0   $, in $ \mathbf{V}_{1}$ the cardinal  $\lambda$ is still
    super-compact and $\boxdot_\lambda$ from \ref{a8} holds. 
    
    (2) Moreover (in part (1)), if (in $ \mathbf{V} _0 $) $h$ is a $\dc $-Laver  diamond  and  $\mathbf{V}_{1}  \models  \lqq\mathbb{P} $ is $\dc$-directed complete'' (see \ref{z41}(3)) \then \, in $ \mathbf{V}_{1}^{\mathbb{P}}$ ($\lambda$ is still  supercompact and) the statement $ \boxdot_{< \lambda, \bar{\theta} _*, \bar{ \theta }}$ holds also  in $\mathbf V_{1}^{\bbP}.$ 
\end{lemma}

\begin{remark}\label{a23}\  

    (1) The following is a major point in \ref{a20} and has caused some confusion. In the definition of $ \boxdot _\lambda $ we  demand only that  the forcing  notion is  $(< \lambda )$-strategically complete. 
    
    To clarify, note:
    
    \begin{enumerate} 
        \item[(A)]  In the proof of $ \boxplus _ \lambda $ holding (i.e. \ref{a20}(1)) we use Easton support iteration  $ \langle  \mathbb{P} _ \varepsilon, \name{ \mathbb{Q} }_ \varepsilon: \varepsilon < \lambda  \rangle$ where $ \mathbb{Q} _ \varepsilon $ is a  $ \mathbb{P} _ \varepsilon $-name of  a  $(< \lambda _ \varepsilon )$-strategically complete forcing notion  from $ {\mathscr H} ( \lambda ), $ 
        where $ \lambda_ \varepsilon = \min \{  \theta : \theta $ regular$ > \varepsilon \ \text{and} \  \le \| \mathbb{P} _ \varepsilon \| \},$
        
        \item[(B)] For \wilog \,  in the the proof of \ref{a29} we use the relevant forcing being $ (< \lambda )$-directed complete but not  so in \ref{a32}). 
    \end{enumerate} 
    
    (2) We may e.g. restrict $\chi$ to be strong limit.
\end{remark}

\begin{PROOF}{\ref{a20}} 
    (1) This is similar to the proof in Laver \cite{Lv78} 
    using Laver's diamond, see Definition \ref{z17}, but we elaborate.  
    
    By Laver \cite{Lv78} \wilog \, there is a Laver diamond $h:\lambda \rightarrow \cH(\lambda)$. Let $E = \{\theta:\theta < \lambda$ is a strong limit cardinal
    and $\alpha < \theta \Rightarrow h(\alpha) \in \cH(\theta)\}$, clearly a club of $\lambda$ and let $\langle \kappa_\varepsilon:\varepsilon <
    \lambda\rangle$ list $\{\theta \in E:\theta$ is strongly  inaccessible$\}$ in increasing order.
    
    We  now define $\mathbf q_\varepsilon$ and $\bar\chi^\varp$  by induction on $\varepsilon \le \lambda$ such that:
    
    \begin{enumerate}
        \item[$(*)_0$] 
        \begin{enumerate}
            \item[(a)]  $\mathbf q_\varepsilon = \langle \bbP_\zeta,\name{\bbQ}_\xi:\zeta \le \varepsilon,\xi < \varepsilon\rangle$ is an Easton support iteration (so $\bbP_\zeta, \name{\bbQ}_\xi$ do not depend on $\varepsilon$), 
            
            \item[(b)]  $\bbP_\zeta \subseteq \cH(\kappa_\zeta)$, 
            
            \item[(c)]  $\bar\chi^\varp = \langle \chi_\zeta:\zeta < \varp\rangle$ where each $\chi_\xi$ is a regular cardinal $\in [\kappa_\xi,\kappa_{\xi +1})$,  
            
            \item[(d)]  $\name{\bbQ}_\xi \in \cH(\chi_{\xi +1})$ is a $\bbP_\xi$-name of a $(<\kappa_\xi)$-strategically  complete forcing notion, 
            
            \item[(e)]   if $h(\xi) = (\name{\bbQ},\chi)$ and the pair $(\name{\bbQ},\chi)$ satisfies the requirements on $(\name{\bbQ}_\xi,\chi_\xi)$ in clauses (c),(d) then $(\name{ \bbQ}_\xi,\chi_\xi) = h(\xi)$.  
        \end{enumerate}
    \end{enumerate}
    
    Concerning clause (b) which says ``$\bbP_\zeta \subseteq \cH(\kappa_\zeta)$", note that for $\zeta$  a limit ordinal letting $\kappa_{< \zeta} = \cup\{\kappa_\xi:\xi < \zeta\}$ we have $\kappa_{< \zeta}$ is strong limit and:
    
    \begin{itemize}
        \item  if $\kappa_{< \zeta}$ is regular, equivalently strongly inaccessible then $\kappa_{< \zeta} = \kappa_\zeta$ and $\bbP_\zeta = \cup\{\bbP_\xi:\xi < \zeta\}$ and so $\bbP_\zeta \subseteq \cup\{\cH(\kappa_\xi):\xi < \zeta\} = \cH(\kappa_{<\zeta}) = \cH(\kappa_\zeta)$, 
        
        \item  if $\kappa_{< \zeta}$ is singular, then $\bbP_\zeta \subseteq  \cH(\kappa^+_{< \zeta}) \subseteq \cH(\kappa_\zeta)$ as $\kappa_\zeta$
        is inaccessible $> \kappa_{< \zeta}$.
    \end{itemize}
    
    Easily we can carry the induction so $\mathbf q_\lambda$ is well defined, $\bbP_\lambda = \cup\{\bbP_\varp:\varp < \lambda\} \subseteq
    \cup\{\cH(\kappa_\varp):\varp < \lambda\} = \cH(\lambda)$ and ``$\xi < \lambda
    \Rightarrow \bbP_\lambda/\bbP_\xi$ is $(< \kappa_\xi)$-strategically complete" hence $\bbP_\lambda/\bbP_\xi$ adds no new sequence of length
    $< \kappa_\xi$ of ordinals.  Clearly it is enough to prove that in $\mathbf V^{\bbP_\lambda}$ we have $\boxdot_\lambda$.
    
    Toward contradiction assume $\chi,\bbP,\cS = \cS_{\chi,\bbP}$ form a counter-example in $\mathbf V^{\bbP_\lambda}$,  hence there are $p_* \in \bbP_\lambda$ and $\bbP_\lambda$-names $\name\chi,\name{\bbP},\name{\cS},\name E$ such
    that $p_* \Vdash_{\bbP_\lambda} ``\name\chi > \lambda$ is regular, $\name{\bbP} \in \cH(\name\chi)$  is $ (< \lambda )$-strategically complete  and $\name{\cS}_{\name\chi,\bbP}$ is defined as in $\boxdot_\lambda$ and $\name E \subseteq [\cH(\chi)^{\mathbf V[\bbP_\lambda]}]^{< \lambda}$ is a club disjoint to $\name{\cS}"$.
    
    As we can increase $p_*$, \wilog \, $\name\chi = \chi$ and let $x = (\chi,\name\bbP)$; and as $\mathbf V \models ``\lambda$ is super-compact and $h$ is a Laver diamond"  for some $(I,D,\mathbf M,\mathbf j,\bfj_0,\bfj_1)$ we have:
    
    \begin{enumerate}
        \item[$(*)_1$]
        \begin{enumerate}
            \item[(a)]  $\mathbf M$ is a transitive class,
            
            \item[(b)]  $\mathbf M$ is a model of ZFC,
            
            \item[(c)]  ${}^\chi \mathbf M \subseteq \mathbf M,$
            
            \item[(d)]  $\mathbf j$ is an elementary embedding from  $\mathbf V$ into $\mathbf M,$
            
            \item[(e)]  $\crit(\mathbf j) = \lambda,$
            
            \item[(f)]  $\mathbf j(h)(\lambda) = (\chi,\name{\bbP}),$
            
            \item[(g)]  $I = [\cH(\chi_1)]^{< \lambda}$ and $\chi_1 > \chi,$
            
            \item[(h)]  $D$ is a fine normal ultrafilter on $I,$
            
            \item[(i)]  $\bfj_0$ is the canonical elementary embedding of $\bfV$ into $\bfV^I/D,$
            
            \item[(j)]  $\bfM$ is the Mostowski Collapse of $\bfV^I/D,$
            
            \item[(k)]  $\bfj_1$ is the canonical isomorphism from $\bfV^I/D$ onto $\bfM,$
            
            \item[(l)]  $\bfj = \bfj_1 \circ \bfj_0$.
        \end{enumerate}
    \end{enumerate}
    
    Moreover, by Definition \ref{z17},
    
    \begin{enumerate}
        \item[$(*)_2$]  $x = \bfj_1(\langle (\sup(u \cap \lambda): u \in I\rangle/D)$.
    \end{enumerate}
    
    Let $\mathbf q = \mathbf j(\mathbf q_\lambda)$ so $\mathbf q = \langle \bbP_\zeta,\name{\bbQ}_\xi:\zeta \le \mathbf j(\lambda),\xi < \mathbf j(\lambda)\rangle$ and $\zeta < \lambda \Rightarrow  \bbP^{\mathbf q}_\zeta = \bbP_\zeta$, etc.
    
    So,
    
    \begin{enumerate}
        \item[$(*)_3$] in $\bfM$ the pair $x = (\chi,\name\bbP)$ satisfies:
        
        \begin{enumerate}
            \item[(a)] $\chi \in (\lambda,\bfj(\lambda))$, 
                  
            \item[(b)] $\name\bbP \in \cH(\chi),$
            
            \item[(c)] $\name\bbP$ is a $\bbP_\lambda$-name of a $(< \lambda)$-strategically complete forcing notion.
        \end{enumerate}
    \end{enumerate}
    
    [Why?  Because $[\bfM]^\chi \subseteq \bfM$ hence $\cH(\chi^+)^{\bfV}
    \subseteq \bfM$].
    
    Now,
    
    \begin{enumerate}
        \item[$(*)_4$]  the following sets belong to $D$:
        
        \begin{enumerate}
            \item[(a)]  $\cS_1 = \{u \in I:x \in u$ and $(\cH(\chi_1),\in) \rest u \prec (\cH(\chi_1),\in)\}$
            
            \item[(b)]  $\cS_2 = \{u \in \cS_1:u \cap \lambda$ is an inaccessible cardinal we call $\lambda_u\},$
            
            \item[(c)]  $\cS_3 = \{u \in \cS_2$: the Mostowski Collapse of $(\cH(\chi_1),\in) \rest u$ is isomorphic, for some $ \chi ^ \dagger$ to $(\cH(\chi^ \dagger ), \in) \}$.  
        \end{enumerate}
    \end{enumerate}
    
    [Why?  As $D$ is a fine and normal ultrafilter on $I$.]
    
    \begin{enumerate}  
        \item[$(*)_5$]  if $ u \in {\mathscr S } _3 \subseteq I $  then let $ \chi _u $ be the $ \chi $ guaranteed to exist by  $ u \in {\mathscr S } _3 $ and  let $ \mathbf{j} _u $ be the Mostowski collapse of $(\cH(\chi_1),\in) \rest u $ onto $ ({\mathscr H} (\chi_u), \in )$
    \end{enumerate} 
    
    \begin{enumerate}
        \item[$(*)_6$]  for every formula $\varphi = \varphi(-) \in  \bbL(\{\in\})$ the following are equivalent:
        
        \begin{enumerate}
            \item[(a)]  $(\cH(\chi_1),  \in) \models \varphi[x]$, 
            
            \item[(b)]  $(\cH(\chi_1),\in)^I/D \models \varphi[\langle h(u \cap \lambda):u \in I\rangle/D]$, 
            
            \item[(c)]  $\cX^1_\varphi \in D$ where $\cX^1_\varphi = \{u \in I:x \in u$ and $(\cH(\chi_1),\in) \rest u \models \varphi[x]\}$, 
            
            \item[(d)]  $\cX^2_\varphi \in D$ where $\cX^2_\varphi = \{u \in I:p_*,x \in u$ and $(\cH(\chi_u),\in) \models \varphi[\bfj_u(x)]\}$ on, 
            
            \item[(e)]  $\cX^3_\varphi \in D$ where $\cX^3_\varphi = \{u \in I:x \in u,\chi_u = \otp(\chi \rest u)$ and $(\cH(\chi_u),\in) \models \varphi[\bfj_u(x)]\}$.
        \end{enumerate}
    \end{enumerate}
    
    [Why? We have (a) $\Leftrightarrow$ (c) as $D$ is a fine normal ultrafilter on $I = \cH(\chi_1)$; we have (c) $\Leftrightarrow$ (d) as $ \mathbf{j}_u$ is an isomorphism from $(\cH(\chi_1),\in) \rest u$ onto  $\cH(\chi_u)$;  we have (d) $\Leftrightarrow$ (e) by the choice of $D$; lastly, (b) $\Leftrightarrow$ (c) by \L os theorem.]
    
    Hence, 
    
    \begin{enumerate}  
        \item[$(*)_7$]  there is $N$ as required in $\bfV^{\bbP }$.  
    \end{enumerate}  
    
    [Why?  Choose $u \in I$ which belongs to all the sets from $D$ mentioned in $(*)_4 + (*)_6$. Let $\zeta = u \cap \lambda$, so it is inaccessible, even measurable, and $\bfj_u(x) = \bfj_u(\chi,\name{\bbP}) = h(\zeta)$
    so (by the choice of $\bfq$) $h(\zeta) = (\chi,\name{\bbQ}_\zeta)$ and $\name{\bbQ}_\zeta$ is a $\bbP_{\bfq,\zeta}$-name.
    
    Let $\bfG$ be a  generic  subset of $\bbP_{\bfq} = \bbP_\lambda$ to which $p_*$ belongs, $\bfG_\zeta = \bfG \cap \bbP_{  \mathbf{q}, \zeta}$, hence it  is a generic subset of $\bbP_{\bfq,\zeta}$ over $\bfV$ hence a generic subset of $\bfj_u(\bbP_{\bfq}) \in \cH(\chi_\zeta)$ and let $N = ((\cH(\chi_1),\in) \rest u)[\bfG],\bbA = (\cH(\chi_\zeta)^{\bfV[\bfG_\zeta]},\in),M =  \bbA^{\name{\bbQ}_\zeta[\bfG_\zeta]}$.  Easily $N$ is as promised, contradiction to the choice of $p_*$.]
    
    So we are done proving part (1).
    
    (2) Let $\bbQ$ be a forcing notion in ${\mathbf V}$ which is $\dc  $-directed complete,  $ \name{ \mathbb{P} } $  is a $ \mathbb{Q} $-name  of a $ (< \lambda )$-strategically complete forcing notion. 
    Let   $\chi_1$  be  large enough so that $\lambda, \mathbb{Q}, \name{\bbP}  \in \cH(\chi_1)$ and it suffices to prove that in ${\mathbf V}^{\bbQ}$, the 
    set $\cS_{\chi_1,\bbQ}$ is stationary. So let $\name E$ be $\bbQ$-name and let  $p \in   \bbP$ be such that  $p \Vdash_{\bbQ}  ``\name E$ a club of $[\cH(\chi_1)]^{< \lambda}$ disjoint to $\name\cS_{\chi_1,\name{\bbP}}$", no need to use a name for $\chi_1$ as we can increase $p$.
    
    Let $\chi \gg \chi_1$; now $\bbQ * \name{\bbP} \in \cH(\chi)$ is a $(< \lambda)$-strategically complete forcing notion and \wilog \,  codes   
    $(\chi_1,E)$.  
    
    As $ \zdc $ holds in $\mathbf V$ we can apply
    it to the forcing $\bbP_{\ge p} * \name{\bbQ}$ so we can find a tuple $(N,\lambda_N,\chi_N,\mathbf j_N,\bbA_N,M_N,\mathbf G_N)$ witnessing it,
    in particular, $(p,\emptyset) \in \mathbf G_N,\bbP * \name{\bbQ} \in N$  so $\chi_1,\name E \in N$.  
    Let $\mathbf G_{\bbP}$ be a subset of $\bbP$ generic
    over $\mathbf V$ which extends  $\{p':(p',\name{ q}') \in \mathbf G_N $ for some $ \name{ q } ' \}$,   
    possible because $\mathbf G_N$ is in $\mathbf V$, a subset of $\bbP$ which has an upper bound, this is the only place we use ``$\bbP$ is  $ \dc $-\dic". 
    
    Next, let $\mathbf V_1 = \mathbf V[\mathbf G_{\bbP}],N_1 = N[\mathbf G_{\bbP}],E_1 = \name E[\mathbf G_{\bbP}],\bbA_1 = \bbA[\mathbf j''_N(\mathbf
    G_{\bbP} \cap N)] = \bbA[\{p':(p',\name q') \in \mathbf G_N  \}],\mathbf G_1 = \{\name q[\mathbf j(\name{ q })
    :(p,\name q) \in \mathbf G_{N}\}.$
    
    Now recalling $p$ forces $\name E$ is disjoint to $\name{\cS}$ clearly, 
    
    \begin{enumerate}
        \item[$(*)$]  $N_1 \in E_1$.  
    \end{enumerate}
    
    Hence,
    
    \begin{enumerate}
        \item[$(*)$]  $N_1 \notin \cS $. 
    \end{enumerate}
    
    But easily in $\mathbf V_1$ we have: $(\lambda_N,\chi_N,\mathbf j_1, \bbA_1,M_1 = M,\mathbf G_1)$ witnesses $N_1 \in \cS \cap E_1$, a contradiction to the choice of $\name E$.
\end{PROOF}

\begin{discussion}\label{a26}
    Suppose that one wishes to force an inequality between two cardinal characteristics. There are two general approaches, which can be labeled as Top-down and Bottom-up.  In the Bottom-up strategy one begins with a universe in which many characteristics are small, e.g. by assuming $2^\lambda = \lambda^+$, and then increases some of them while trying to keep the smallness of the rest.  In the Top-down strategy one begins with a universe in which many characteristics are large.  The forcing aims to decrease some of them while keeping the large value of the rest.
    
    We shall use the Top-down approach, so we begin by increasing $\gb_\lambda$ (and $\gd_\lambda$) to some $\mu = \cf(\mu) > \lambda$. Notice that $\gb_\lambda$ is a relatively small characteristic and,
    in particular, always $\gb_\lambda \le \gd_\lambda$.  The next step will be to decrease $\cov_\lambda(\meagre)$ in such a way that maintains the fact that $\gd_\lambda = \mu$.  We shall increase
    $\gb_\lambda$ by using the generalization to $\lambda$ of Hechler forcing.  This is a standard way to achieve this goal, but we spell out the proof since it demonstrates the way that we employ Lemma \ref{a20}. 
\end{discussion}

\begin{claim}\label{a29}
    Assume that:
    \begin{enumerate}
        \item[(a)]  $\lambda$ is super-compact,
        
        \item[(b)]  $\lambda < \mu = \cf(\mu) = \mu^\lambda$ and,
        
        \item[(c)] $ \zdc $  and $ \theta_* >  (2^\mu )^+.$
    \end{enumerate}
    
    \Then \, one can force $\gb_\lambda = \gd_\lambda = \mu$ while keeping the super-compactness of $\lambda$ and the principle $\zdc$ and \ref{a20} still holds by \ref{z47}(2). 
\end{claim}

\begin{PROOF}{\ref{a29}}
    Begin with the preparatory forcing of Lemma \ref{a20} to make $\lambda$ indestructible  and to force $\zdc $, hence it will be preserved by any further $(< \lambda)$-\dic forcing.  By \ref{a20} as in the applications of Laver-indestructibility we can assume that $\GCH$ holds above $\lambda$ after the preparatory forcing.  In particular, if $\mu = \cf(\mu) > \lambda$ then $\mu^\lambda = \mu$ follows.
    
    Let $\bbD_\lambda$ be the generalized Hechler forcing.  So a condition $p \in \bbD_\lambda$ is a pair $(\eta_p,f_p)$ such that $\eta_p \in {}^{<\lambda} \lambda,f_p \in {}^\lambda \lambda$ and $\eta_p \trianglelefteq f_p$.  If $p,q \in \bbD_\lambda$ then $p \le q$ iff $\eta_p \trianglelefteq \eta_q$ and $f_p(\alpha) \le f_q(\alpha)$ for every $\alpha \in \lambda$.
    
    Let $\bfq = \langle \bbP_\alpha,\name{ \bbQ}_\beta:\alpha \le \mu,\beta < \mu \rangle$ be a $(<\lambda)$-support iteration of the generalized Hechler forcing notions for $\lambda$.  Explicitly, $\name{\bbQ}_\alpha$ is the $\bbP_\alpha$-name of $\bbD_\lambda$ in $\bfV^{\bbP_\alpha}$ for every $\alpha < \mu$.  Denote the generic $\lambda$-Hechler for $\name{\bbQ}_\alpha$ by $\name f^*_\alpha$.  So $\bbP_\mu$ is the limit and choose a generic $\bfG \subseteq \bbP_\mu$.  We claim that
    $\bfV[\bfG] \models ``\gb_\lambda = \gd_\lambda = \mu"$ as witnessed by $\langle f^*_\alpha:\alpha < \mu\rangle$.  Notice that $2^\lambda = \mu$ in $\bfV[\bfG]$, so it is sufficient to prove that $\gb_\lambda = \mu$ in $\bfV[\bfG]$.
    
    Since $\lambda$ is regular, each $\name{\bbQ}_\alpha$ is $(< \lambda)$-complete  and even $ (< \lambda )$-directed complete.  By Fact \ref{z47}(2)  $\bbP_\alpha$ is $(<\lambda)$-directed  complete as well, for every $\alpha \le \mu$.  Likewise, each $\name{\bbQ}_\alpha$ satisfies $*^\omega_\lambda$  so $\bbP_\mu$ is $\lambda^+$-c.c. (see \cite{Sh:80} or \cite{Sh:1036}).  It follows that $\bfV[\bfG]$ preserves cardinals and cofinalities.  Moreover, no new $(< \lambda)$-sequences of ordinals are introduced.  Notice also that $\bbP_\mu$ is $(< \lambda)$-\dic and so   (by \ref{z47}(3)) it is $\dc $-\dic. 
    
    By \ref{a20}(2) this implies $\bfV[\bfG] \models ``\lambda$ is super-compact and $\zdc$ holds".
    
    The main point is that $\{f^*_\alpha:\alpha < \mu\}$ is a cofinal family in $({}^\lambda \lambda)^{\bfV[\bfG]}$.  For this, assume that $\Vdash_{\bbP_\mu} ``\name f \in {}^\lambda \lambda"$. For every $\alpha < \lambda$ fix a maximal antichain  $\langle p_{\alpha, i } :i < i_\alpha \le \lambda\rangle$ of conditions which force a value to $\name f(\alpha)$.  Let $\delta = \sup(\cup\
    \{\dom(p_{\alpha , i}):\alpha < \lambda,i < i_\alpha \})$. 
    
    Since $\lambda < \mu = \cf(\mu)$ we see that $\delta < \mu$, and clearly $\name f$ is a $\bbP_\delta$-name.  We conclude, therefore, that $\name f$ is dominated by $\name f^*_{\delta +1}$ and hence $\{f^*_\alpha:\alpha < \mu\}$ exemplifies $\gb_\lambda = \mu$.  This fact completes the proof.
\end{PROOF}

Our second lemma is the main burden of the proof.  The statement of the theorem requires $\lambda$ to be super-compact, in order to obtain the indestructibility properties given by Lemma \ref{a20}.  The
combinatorial part given in Lemma \ref{a32} below requires only strong inaccessibility.  However, we assume super-compactness in order to keep  $ \zdc $,  

\begin{lemma}\label{a32}
    Assume that:
    
    \begin{enumerate}
        \item[(a)]  $\lambda$ is super-compact,  
        
        \item[(b)] $\lambda < \mu = \cf(\mu) = \mu^{\lambda},$  
        
        \item[(c)]  $\zdc $ and $ \theta _* \ge (2^\mu )^+$.
    \end{enumerate}
    
    \Then \, there exists a $\lambda^+$-c.c. $(< \lambda)$-strategically  forcing notion $\bbP$ such that $\Vdash_{\bbP} ``\gd_\lambda = \mu \wedge \cov_\lambda(\meagre) = \kappa";$ also $\Vdash_{\bbP}$``$\lambda$ is supercompact''. 
\end{lemma}

\begin{PROOF}{\ref{a32}}
    By claim \ref{a29} \wilog \, $\gb_\lambda = \gd_\lambda = \mu$.  In particular, $\lambda$ is super-compact and $\boxdot_{<\lambda, \theta_{*}, \overline{\theta}}$ holds in the generic extension.  Let $\langle f^*_\alpha:\alpha <\mu\rangle$ witness $\gb_\lambda = \gd_\lambda = \mu$ and \wilog \, $\alpha < \beta < \mu \Rightarrow f^*_\alpha <_{J^{\bd}_\lambda} f^*_\beta$.
    
    Recalling Definitions  \ref{z32}, \ref{z33}, Claim \ref{z35}, Theorem \ref{z38},
    in $\mathbf V$ there are $\beta(*), \mathbf q,\bar u,\cU_*,\ldots$ such that:
    
    \begin{enumerate}
    \item[$(*)_1(A)$]  $\mathbf q \in \mathbf Q_{\lambda,\bar\theta,\beta(*)}$ so in
      particular we have (in $\bfq$):
    
    \begin{enumerate}
        \item[(a)]  $\langle \bbP_{0,\alpha},\name{\bbQ}_{0,\beta}:
        \alpha \le \beta(*),\beta <
        \beta(*)\rangle$ is a $(< \lambda)$-support iteration,  
        
        \item[(b)]  $\bar u = \langle u_\beta:\beta <
        \beta(*)\rangle,\bar{\cP} = \langle \cP_\beta:\beta < \beta(*)\rangle$, 
        
        \item[(c)]  $u_\beta \subseteq \beta,\cP_\beta \subseteq [u_\beta]^{\le \lambda}$ is closed under subsets, 
        
        \item[(d)]  $\name{\bbQ}_{0,\beta}$ has generic 
        $\name\eta_\beta \in \prod\limits_{\varepsilon < \lambda} \theta_\varepsilon$, and is $(< \lambda)$-strategically complete, 
        
        \item[(e)]  $\name\bbQ_{0,\beta}$ is as in \ref{z32}(e) so is $\subseteq \bbQ_{\bar\theta}^{\mathbf V[\langle \name \eta_\alpha:\alpha \in u_\beta\rangle]}$ and $\Vdash_{\bbP_{\beta +1}} ``\name\eta_\beta \in \prod\limits_{\varepsilon < \lambda} \theta_\varepsilon"$ and $\name{\bar\eta} = \langle \name\eta_\beta:\beta < \beta(*)\rangle$, 
        \item[(f)]  $\cU_* \subseteq \beta(*)$ has order type $ \ggk $ and $\langle \beta^*_i:i \le \kappa\rangle$ lists $\cU_* \cup \{\beta(*)\}$ in increasing order, and $\beta(*) = \sup(\cU_{*}),$ 
        \item[(g)]   if $\beta \in \cU_*$ then    $\cU_* \cap  \beta  \subseteq u_ \beta$ and $ [\cU_ * \cap \beta ]^{\le \lambda} \subseteq \cP_\beta$  and $\Vdash_{\bbP_{0,\beta +1}}$ ``if $\nu \in \mathbf V[\langle \name\eta_\alpha: \alpha \in \cU_* \cap \beta\rangle] \cap  \prod\limits_{\varepsilon < \lambda} \theta_\varepsilon$ then $\nu <_{J^{\bd}_\lambda} \name\eta_\beta"$, 
        
        \item[(h)]  if $\alpha \le \beta(*)$ then $\bbP_{0,\alpha}$ is $(< \lambda)$-strategically complete, $ \dc $-\dic and $\lambda^+$-c.c.,  
        
        \item[(i)] ${\bbP}_{1, \alpha}, {\bbP}_{1, \cU}$ are as in \ref{z33},   
        
        \item[(j)]  $ \theta _*  \ge (2^\mu  )^+$.  
    \end{enumerate} 
    
    \item[(B)]  letting $\bbP'_i = \bbP_{\mathbf q,1,\{\beta^*_j:j <i\}}$ for $i \le \kappa   $  
    we have: 
    
    \begin{enumerate}
        \item[(a)]  The sequence $ \bar{ \mathbb{P} }'= \langle \bbP'_i:i \le \ggk\rangle$ of forcing notions is $\lessdot$-increasing, and is continuous for 
        ordinals $i \le \ggk$ of cofinality $ > \lambda $ see \cite[2.5(8)=Lb14(8)]{Sh:1126}, but the continuity will not be used, 
        
        \item[(b)]  $\bbP'_i$ is $(< \lambda)$-strategically complete for $i \le \ggk$, 
        
        \item[(c)]  $(\prod\limits_{\varepsilon < \lambda} \theta_\varepsilon)^{\mathbf V[\bbP'_{\ggk}]} = \cup\{(\prod\limits_{\varepsilon < \lambda} \theta_\varepsilon)^{\mathbf V[\bbP'_i]}:i < \ggk\}$, 
        
        \item[(d)]  The sequence $\langle \bbP_{1,\beta}:\beta \le \beta(*)\rangle$ is a sequence of forcing notions, $\lessdot$-increasing and if $\beta \le \beta(*)$ then $\bbP_{0,\beta} \lessdot \bbP_{1,\beta}$, in fact is dense in it and if $i \le \ggk$ then $\bbP'_i \lessdot \bbP_{1,\beta^*_i}$.
    \end{enumerate}
    \end{enumerate}
    
    We shall mention more properties later.
    
    [Why are there such objects?  We apply \ref{z38} and
    \ref{z32} and \ref{z35}, that is \cite{Sh:1126}].   
    
    Also,
    
    \begin{enumerate}
        \item[$(*)_2$] 
        \begin{enumerate}
            \item[(a)] recall $\langle \beta^*_i:i \le \ggk  \rangle$ lists $\cU_* \cup \{\beta(*)\}$ in increasing order, 
            
            \item[(b)]  for $i < \ggk $ let $\name g'_i$ be $\name\eta_{\beta^*_i}$ (to avoid excessive subscripts), 
            
            \item[(c)]   let $\name{\bar g}' = \langle \name g'_i:i < \kappa\rangle$, 
            
            \item[(d)]  $\cP_\alpha = \cP_{\mathbf q,\alpha}$ and \wilog \, $u_\alpha = \cup\{u:u \in \cP_\alpha\}$ for $\alpha < \beta(*)$.
        \end{enumerate}
        
        \item[$(*)_3$]   if $u \in \cP_\alpha,\alpha < \beta(*)$ then $\Vdash_{\bbP_{0,\alpha +1}} ``\name \ggy _\alpha \in \prod\limits_{\varepsilon < \lambda} \theta_\varepsilon$ dominates $(\prod\limits_{\varepsilon < \lambda} \theta_\varepsilon)^{\mathbf V[\langle \name g_\beta:\beta \in u \rangle]}"$, the order being modulo $J^{\text{bd}}_\lambda$.
    \end{enumerate}
    
    [Why? By the choice of the forcing, see \ref{z23}
    or $(*)_1(A)(g)$ above].
    
    \begin{enumerate}
        \item[$(*)_4$]   we have
        
        \begin{enumerate} 
            \item[(a)]  $\Vdash_{\bbP'_\kappa} ``\name{\bar g}' = 
            \langle \name g'_i:i < \kappa\rangle$ is 
            $<_{J^{\text{bd}}_\lambda}$-increasing and
            cofinal in $(\prod_{\varepsilon < \lambda} \theta_\varepsilon, <_{J^{\text{bd}}_\lambda})"$.
            
            \item[(b)] for every $p \in \bbP_{\beta(*)}$ and $\alpha \in \dom(p) \cap \cU_{*},$ for every large enough $i < \kappa,$ we have $p \Vdash_{\bbP_{0, \beta(\ast)}}$``$\name{f}_{p(\alpha)} \leq \name{g}_{i}' = \name{\eta}_{\beta_{i}^{*}} \mod J_{\lambda}^{\rm{bd}}$'',
            
            \item[(c)] $\Vdash_{\bbP_{0, \beta(*)}}$``$ \overline{\name{g}}' = \langle \name{g}_{i}': i < \kappa \rangle$ is $<_{J_{\lambda}^{\rm{bd}}}$-increasing and cofinal in $\left( \Pi_{\varp < \lambda \theta_{\varp}} , <_{J_{\lambda}^{\rm{bd}}}\right)$''.  
        \end{enumerate} 
    \end{enumerate}
    
    [Why? Clause (a) holds by  $(*)_3$ noting that $(\prod\limits_{\varepsilon < \lambda}  \theta_\varepsilon)^{\mathbf V[\bbP'_\kappa]} = 
    \cup\{(\prod\limits_{\varepsilon < \lambda} 
    \theta_\varepsilon)^{\mathbf V[\bbP'_i]}:i < \kappa\}$ which holds by   \ref{z38}(A)(d).  Clause (b) holds  by \ref{z38}(A)(e) alternatively, in \cite{Sh:1126} this notation means:
    
    \begin{itemize}
        \item  if $ \alpha \in \cU_{*}$ and $t \in L_{\bfm} \setminus M_{\bfm}, \alpha(i) \in u_{\alpha} \cap (t / E_{\bfn}')$ for $i < \lambda$ and $\bfB$ a $\lambda$-Borel function from ${}^{\lambda} \left( \Pi_{\varp < \lambda} \theta_{\varp} \right)$ into $\Pi_{\varp < \lambda} \theta_{\varp}$ \underline{then} for every $i < \kappa$ large enough $\Vdash$``$\bfB(\dots, \name{\eta}_{\alpha(\varp)}, \dots)_{i < \lambda} \leq_{J_{\lambda}^{\rm{bd}}} g'$''. 
    \end{itemize}.
    
    Clause  (c) holds by $(*)_{1}$(A)(g)].  
    
    Now, 
    
    \begin{enumerate}
        \item[$(*)_5$]   $\Vdash_{\bbP'_\kappa}$
        ``cov$_\lambda$(meagre) $\le \kappa$".
    \end{enumerate}
    
    [Why?  First, notice that we can look at $\prod\limits_{\varepsilon < \lambda} \theta_\varepsilon$ instead of ${}^\lambda 2$.  
    
    Second, for each $\varepsilon <  \lambda,i < \kappa$ the set $B_{\varepsilon,i} = \{\eta \in \prod\limits_{\xi < \lambda} \theta_\xi$: for every $\zeta \in [\varepsilon,\lambda)$ we have $\eta(\zeta) \le \name g'_i(\zeta)  < \theta_\zeta\}$ is closed nowhere dense, and by $(*)_4$ we have $\mathbf V^{\bbP'_\kappa} \models ``\prod_{\zeta < \lambda} \theta_\zeta = \cup\{B_{\varepsilon,i}:\varepsilon < \lambda,i < \kappa\}"$.  In fact, $\langle B_{0,i}:i < \kappa\rangle$ suffice.
    
    Alternatively we have $\langle g'_i : i < \kappa \rangle $ is $ < _{J^{\rm bd}_\lambda  } $-increasing cofinal in $ \Pi _{\varepsilon < \lambda } \theta _ \varepsilon $ and let $ {\cW }   _{i, \zeta } := \{ \eta: \eta \in {}^ \lambda 2 $ and for every $ \varepsilon \in [ \zeta, \lambda )$ we have either $\eta \restriction [ \Sigma _{\xi  < \varepsilon} \theta _  \xi , \Sigma _{\xi  \le  \varepsilon }  \theta _  \xi)$ is constantly zero
    or ${\rm min } \{ \alpha : \Sigma _{\xi  < \varepsilon }  \theta _  \xi + \alpha \in \eta ^{-1}(\{ 1 \} ) \} < g'_i (\varepsilon ) \}.$ So $ {\cW }   _{i, \zeta } $ is a closed nowhere dense subset of $ {}^ \lambda  2 $ and $ \cup \{ {\cW }   _{i, \zeta } : i  <  \kappa, \zeta < \lambda \} = {}^\lambda 2 $ and $ \kappa \times \lambda$ has cardinality $ \lambda + \kappa = \kappa $   because if $ f \in {}^{ \lambda } 2  $ then we define  $\nu _f  \in \Pi  _{\varepsilon < \lambda }  \theta _ \varepsilon $  as follows:  for $ \varepsilon < \lambda$:
    
    \begin{enumerate} 
        \item[(a)]  if $ f \rest [ \Sigma _{\xi < \varepsilon }\theta _ \xi, \Sigma _{\xi \le \varepsilon } \theta _ \varepsilon )$ is not  constantly zero then  we let $ \nu _f (\varepsilon) = \min \{\alpha : f(\Sigma _{\xi < \varepsilon } \theta _ \varepsilon  + \alpha ) = 1 \}$;
                
        \item[(b)] if otherwise then  let $ \nu_f (\varepsilon ) = 0.$
    \end{enumerate} 
    
    So there are $ i < \kappa $ and $ \varepsilon < \lambda $ such that: $ \zeta \in [\varepsilon, \lambda ) \Rightarrow \nu _ f ( \zeta ) < g'_i(\zeta ) $. Now it is  easy to check that $ f \in {\mathscr W } _{i, \varepsilon}.$]  
    
    Lastly, 
    
    \begin{enumerate}
        \item[$(*)_7$] $\Vdash_{\bbP'_\kappa}``\cov_\lambda(\meagre) \ge
        \kappa"$.
    \end{enumerate}  
    
    [Why?  For $ i < \kappa$ let us define the $\bbP'_{i +1}$-name $\name\nu '_i$ of a member of ${}^\lambda 2$ by $\name \nu '_i(\varepsilon)=0$  iff $\name g'_i(\varepsilon)$ is even.  Now  clearly $\Vdash_{\bbP'_{i +1}} ``\name   \nu '_i$ is a  $\lambda$-Cohen sequence over $\mathbf V^{\bbP'_i}"$. (But let us  elaborate; $\name  \nu  '_i$ is also a $\bbP_{\beta^*_i+1}$-name and $\Vdash_{\bbP_{\beta^*_i+1}} ``\name \nu '_i$ is $\lambda$-Cohen over $\mathbf V^{\bbP_{\beta^*_i}}$ hence over $\mathbf V^{\bbP'_i}"$; the last hence because $\bbP'_i \lessdot \bbP_{1,  \beta^*_i}$.  As $\bbP_{\beta^*_i+1} \lessdot \bbP_{\beta^*_{i+1}}$ and $\bbP'_{i+1} \lessdot \bbP_{\beta^*_{i+1}}$ we are done 
    {)}.  
    
    Also every closed nowhere dense subset of ${}^\lambda 2$ from $\mathbf V^{\bbP'_{\ggk}}$ is from $\mathbf V^{\bbP'_i}$ for some $i < \ggk$.   So if $p \Vdash ``\cov_\lambda(\meagre) < \kappa"$ then for some $\zeta < \kappa$ and $\name A_\varp(\varp < \zeta)$ we have $p \Vdash ``\name A_\varp$ is a closed no-where dense subset of ${}^\lambda 2$ for $ \varepsilon < \zeta"$   and $p \Vdash ``\bigcup\limits_{ \varepsilon < \zeta } \name{ A}_ \varepsilon$ is equal to the set of ${}^\lambda 2"$.  \Wilog \, each $\name A_\varp$ is a $\bbP_{i(\varp)}$-name, $i(\varp) < \kappa$
    and recall that $\kappa$ is regular.
    Hence $i = \sup\{i(\varp):\varp < \zeta\} < \kappa$ and $\name g '_i$ gives a  contradiction to the choice of $\langle \name A_\varp:\varp <
    \zeta\rangle$; so $(*)_{6}$ holds 
    indeed.]
    
    The reader may look at some explanation in \ref{a41}.
    
    Now we come to the main and last point recalling $\langle f^*_\alpha:\alpha < \mu\rangle$ from Claim \ref{a29}
    
    \begin{enumerate}
        \item[$(*)_7$]  $\Vdash_{\bbP'_{\ggk}}$ ``no 
        $\name f \in ({}^\lambda \lambda)$ dominates 
        $\{f^*_\alpha:\alpha < \mu\}"$.
    \end{enumerate}
    
    We shall show that it suffices to prove $(*)_7$ for proving Lemma  \ref{a20}(2), and  then   that $(*)_7$ holds, thus finishing.
    
    Why it suffices?  As $\langle f^*_\alpha:\alpha < \mu\rangle$ is $<_{J^{\text{bd}}_\lambda}$-increasing and cf$(\mu) = \mu > \lambda$, this implies $\Vdash_{\bbP'_\kappa} ``\gd_\lambda \ge \mu"$.  Also in $\mathbf V,\mu^\lambda = \mu > \kappa > \lambda$ and $|\bbP'_{\ggk}| = \kappa^\lambda$ by (A)(g) of \ref{z35}(4) which is $\le \mu$ and $\bbP'_\kappa$ satisfies the $\lambda^+$-c.c. hence $\Vdash_{\bbP'_\kappa} ``2^\lambda = \mu  " $,   hence together $\Vdash_{\bbP'_\kappa} ``\gd_\lambda = \mu"$.   Also by $(*)_1(B)(b)$, ``$\bbP'_{\ggk}$ is $(< \lambda)$-strategically complete  and  $ \lambda^+$-c.c."  and by $(*)_5 + (*)_6$ we know that  ``$\cov_\lambda(\meagre) = \kappa"$ so we are
    done; hence $(*)_7$ is really the last piece missing.  
    The rest of the proof is dedicated to proving that $(*)_7$ holds.  
    
    We shall use further nice properties of $\bbP'_j,\name g'_i (j \le \ggk,i < \ggk)$ which hold by $(*)_1 + (*)_2$ (and $(*)_3,(*)_4$)  and their proof, i.e. \ref{z35}, \ref{z38} 
    and see \cite[2.12=Lb35, 2.13=Lb38]{Sh:1126}.  
    
    \begin{enumerate}
        \item[$\boxplus_1$] 
        
        \begin{enumerate}
            \item[(a)]
            
            \begin{enumerate}
                \item[$(\alpha)$]  $\langle \name g'_\gamma:
                \gamma < \ggk\rangle$ is generic for $\bbP'_{\ggk}$,  i.e., if $\mathbf G$ is a subset of $\bbP'_{\ggk}$ generic over $\mathbf V$ and $g'_i = \name g'_i[\mathbf G]$ \then \, $\mathbf V[\mathbf G] = \mathbf V[\langle g'_i:i < \ggk\rangle],$
                
                \item[$(\beta)$]  if in addition $\nu \in  ({}^\lambda \lambda)^{\mathbf V[\mathbf G]}$ \then \, for some  $\rho \in ({}^\lambda \ggk)^{\mathbf V}$ and $\lambda$-Borel  function $\mathbf B \in \mathbf V$ we have $\nu =  \mathbf B(\langle g'_{\rho(\varepsilon)}:\varepsilon  < \lambda\rangle)$
            \end{enumerate}
            
            \item[(b)]  if in $\mathbf V[\mathbf G], ''_\gamma \in \prod\limits_{\zeta < \lambda} \theta_\zeta$ for $\gamma < \ggk$ and the set $\{(\gamma,\zeta):\gamma < \ggk$ and $\zeta < \lambda$  and $g''_\gamma(\zeta) \ne g'_\gamma(\zeta)\}$  has cardinality $< \lambda$ then  $\bar g'' = \langle g''_\gamma:\gamma < \ggk\rangle$ is generic for $\bbP'_{\ggk}$ and $\mathbf V[\bar g''] = \mathbf V[\bar g']$;  similarly for $ \mathbb{P} _ {\beta (*)},$  
                
            \item[(b)$^+$]  similarly for any $ \gamma _ \bullet \le  \kappa $  and  $ \langle g''_ \gamma : \gamma < \gamma _ \bullet \rangle $,  really follows, 
            
            \item[(c)]  $\Vdash_{\bbP'_{\ggk}} ``\name g'_\gamma$ dominates $(\prod\limits_{\varepsilon < \lambda} \theta_\varepsilon)^{\mathbf V[\bbP'_\gamma]}"$  for $ \gamma < \ggk $, 
            
            \item[(d)] if $\langle \zeta(\gamma):\gamma < \ggk\rangle$ is an increasing sequence of ordinals $< \ggk$ (from $\mathbf V!$), then $\langle g'_{\zeta(\gamma)}:\gamma < \ggk\rangle$ is generic for $\bbP'_{\ggk}$ (over $\bfV$); 
            
            \item[(e)]  if $\gamma \le \ggk$ then  $\bbP'_\gamma$ is $(< \lambda)$-strategically complete and satisfies the $\lambda^+$-c.c..
        \end{enumerate}
    \end{enumerate}
    
    We shall use $\boxplus_1$ freely, this (mainly $ \boxplus _1(d)$) had been  the  motivation  for \cite{Sh:1126}. 
    
    To prove $(*)_7$ assume toward contradiction that this fails, and hence for some condition $p^* \in \bbP'_{\ggk}$  and $\bbP'_{\ggk}$-name $\name f$ and $\lambda$-Borel function $\mathbf B$  (from $ \mathbf{V} $)  and $\rho_ \bullet  \in {}^\lambda \ggk$ we have:
    
    \begin{enumerate}
        \item[$\circledast_0$]  $p^* \Vdash_{\bbP'_{\ggk}}
        ``\name f \in {}^\lambda \lambda$ and it does  dominate  $ \{ f^*_ \alpha : \alpha < \mu \} $, equivalently  $({}^\lambda \lambda)^{\mathbf V}"$ and $ \name f = \mathbf B(\langle \name g'_{\rho_\bullet (i)}:i < \lambda\rangle)$.
    \end{enumerate}
    
    Now let $\chi$ be regular large enough and we choose $\bar N = \langle N_\varepsilon:\varepsilon < \lambda \rangle$ such that:   
    
    \begin{enumerate}
        \item[$\circledast_1$]  
        \begin{enumerate}
            \item[(a)]  $N_\varepsilon$ is as in $ \zdc $ for the forcing notion $\mathbb{P} _{0, \beta (*)}$  (equivalently $ \mathbb{P} _{1, \beta (*)} $,   \underline{not} $ \mathbb{P} '_\kappa $),  that is $N_\varepsilon \in \cS_{\chi, \mathbb{P} _{1, \beta (*) }}$  see $ \zdc$ of \ref{a8}(2),   
             
            \item[(b)] 
            
            \begin{itemize}
                \item $\lambda_\varepsilon  = \lambda (\varepsilon ) :=  {\rm otp}(N_\varepsilon \cap \lambda ) > \lambda^-_\varepsilon := \Sigma\{\|N_\zeta\|:\zeta < \varepsilon\} \ge \Sigma\{\lambda_\zeta:\zeta < \varepsilon\} \ge \varepsilon $,  
                
                \item $\bar N \restriction \varepsilon \in N_\varepsilon$
                and$ {\rm otp } ( N_ \varepsilon \cap \kappa ) < \theta _ {\lambda ( \varepsilon )},$ moreover $ \otp(\theta _* \cap N_ \varepsilon) < \theta _{\lambda (\varepsilon )} $, (really $ \mathbb{P} _{1, \beta (*) } \in {\mathscr H} ( \theta _* ) $) hence $\bigcup\limits_{\zeta < \varepsilon} N_\zeta \subseteq N_\varepsilon$, 
                
                \item $\lambda_{\varp}$ is an inaccesible cardinal. 
            \end{itemize}
            
            \item[(c)]  $\lambda, \kappa, \mu ,\theta _*, \bar\theta,\mathbf q, \bar{ \mathbb{P} }' , \cU_*, p^*,\name f, \mathbf B,\rho _ \bullet, \bar{ g } ' $ belong to $N_\varepsilon$.
        \end{enumerate}  
    \end{enumerate}
    
    Next choose $f^* \in {}^\lambda \lambda$, i.e. $\in ({}^\lambda \lambda)^{\mathbf V}$, such that:
    
    \begin{enumerate}
        \item[$\circledast_2$] for arbitrarily large $\varepsilon < \lambda$ for some $\zeta \in [\lambda^-_\varepsilon,\lambda_\varepsilon)$ we have $f^*(\zeta) > \lambda_\varepsilon$, (we can demand more: for every large enough $\varepsilon < \lambda $).
    \end{enumerate}
    
    For $\varepsilon < \lambda$ let $(\lambda_\varepsilon,\chi_\varepsilon,\mathbf j_\varepsilon, M_\varepsilon,\bbA_\varepsilon,\mathbf G  ^+  _\varepsilon)$ be a witness  for $(N_\varepsilon,\mathbb{P} _{1, \beta (*)} ,  \chi)$ recalling $\chi$ was chosen after $\circledast_{1}$(a) we have  $\zdc$  from Definition \ref{a8}(2) so $\lambda_\varepsilon \in (\varepsilon,\lambda)$ is strongly inaccessible and $\varepsilon < \zeta < \lambda \Rightarrow \lambda_\varepsilon < \lambda^-_\zeta < \lambda_\zeta$,  recalling $\circledast_1$ and noting 
    $\langle \lambda^-_\varp:\varp <\lambda\rangle$ is an increasing  continuous sequence of cardinals
    below $\lambda$.  Let  $\mathbf{G} ^ \dagger _ \varepsilon = \mathbf{G} ^+_ \varepsilon  \cap \mathbf{j} _\varepsilon ( \mathbb{P} '_{\gamma ( * )})$ recalling $ \mathbf{G} ^+_ \varp \subseteq \mathbf{j} _ \varepsilon (\mathbb{P} _{1, \beta (*)})$  and noting $ \mathbf{G} ^\pigyon _ \varepsilon \subseteq \mathbf{j} _ \varepsilon (\mathbb{P} ' _\ggk) \subseteq \mathbb{A} _ \varepsilon $. 
    
    Let  (for $ \varepsilon < \lambda$):   
    
    \begin{enumerate}
        \item[$\circledast_3$]
            
        \begin{enumerate}
            \item[(a)]  $v_\varepsilon = N_\varepsilon \cap \ggk$, 
            
            \item[(b)]  $\exi_\varp =  \exi(\varepsilon) = \otp(v_\varepsilon)$ and so $\exi_\varepsilon = \mathbf j_\varepsilon(\ggk) $  
            etc, 
            
            \item[(c)]  $\bar\gamma^\varepsilon = 
            \langle \gamma_i(\varepsilon):i < \exi_\varepsilon \rangle$ list $v_\varepsilon$ in increasing order,    
            
            \item[(d)]   for $i < \otp(v_\varepsilon)$, 
            equivalently $i < \mathbf j_\varepsilon(\ggk) = \exi_ \varepsilon$ let $\eta^\varepsilon_i = (\mathbf j_\varepsilon(\name g'_{\gamma_i(\varp)})) [\mathbf{G} ^ \pigyon _ \varepsilon ]  \in  \prod\limits_{\zeta < \lambda_\varepsilon}\theta_\zeta \cap \mathbb{A} [\mathbf{G} ^\pigyon_\varepsilon ] $ ,  
             
            \item[(e)]  let $\bar \eta^\varepsilon = \langle \eta^\varepsilon_i:i < \exi_\varp\rangle$.
        \end{enumerate}
    \end{enumerate}
    
    Note that clearly,
    
    \begin{enumerate}
        \item[$\circledast_4$]   for each $\varepsilon < \lambda$ we have:  
        
        \begin{enumerate}
            \item[(a)]  $\bar \eta^\varepsilon$ is generic
            for $(\bbA_\varepsilon,\mathbf j_\varepsilon(\bbP'_{\ggk}))$, moreover
            
            \item[(b)]  if we change $\eta^\varepsilon_i(\zeta)$ (legally, i.e. to an ordinal $< \theta_\zeta$) for
            $< \lambda_\varepsilon$ pairs $(i,\zeta) \in \otp(v_\varepsilon) \times \lambda_\varepsilon$ and get $\bar \eta'$, \underline{then} also $\bar \eta'$ is generic for $(\bbA_\varp,\mathbf j_\varepsilon (\bbP'_{\ggk}))$, clearly $  N_ \varepsilon [\bar{ \eta } ^ \varepsilon ]  = N_ \varepsilon [\bar{ \eta } ' ]$, 
            
            \item[(c)] there is a unique  $ \mathbf{G}'_ \varepsilon $ a subset of $ \mathbb{P} _{1, \beta (*)}  \cap N_ \varepsilon$ generic over $ N_ \varepsilon$ such that  $\mathbf{j} '' _ \varepsilon (\mathbf{G}' _ \varepsilon )= \mathbf{G} ^+_ \varepsilon $  so $ {\mathscr H} ( \chi _ \varepsilon ) = \mathbb{A} _\varepsilon [\mathbf{j} _ \varepsilon '' (\mathbf{G}'_ \varepsilon)]$
        \end{enumerate}
                
        [Why this equality holds? By $\circledast_{1}$(a) and the first line after $\circledast_{2}$ i.e. \ref{a8}(2); recall  $\mathbf{G} ^ \dagger_\varepsilon = \mathbf{G} ^+_ \varepsilon \cap \mathbf{j}_\varepsilon (\mathbb{P} ' _{\ggk })$  and we have  $(\mathbf{j} _ \varepsilon (\name{ \eta} _{\gamma _i(\varepsilon)}[\mathbf{G}'  _\varepsilon ])= \eta ^ \varepsilon _i$,] 
        
        \begin{enumerate} 
            \item[(d)] like $\boxplus_1$ with $\mathbf V,    \mathbb{P} _{1,\beta (*)}, \mathbf{G},\lambda$ there standing for $\bbA_\varp,\mathbf j_\varepsilon(\mathbb{P} _{1,\beta (*)}) \mathbf{G}^+  _ \varepsilon ,\lambda_\varepsilon$ here.  
        \end{enumerate}
    \end{enumerate}
    
    Hence we have, 
    
    \begin{enumerate}
        \item[$\circledast'_4$]   for $\varepsilon < \lambda$, 
        
        \begin{enumerate}
            \item[(a)]
            
            \begin{enumerate} 
                \item[$( \alpha ) $] let $ \Xi ^ \dagger _ \varepsilon : = \{ \bar{ \nu }: \bar{ \nu } = \langle \nu _i : i < \kappa _\varepsilon  \rangle $ and for some  $ {\mathbf G} ^\dagger  \subseteq {\mathbb   P}'_{ } \cap N_\varepsilon $ generic over $ N_ \varepsilon $
                we have $ \nu _i \in \Pi _{ \xi < \lambda _ \varepsilon} \theta  _  \xi$ satisfies $ \xi < \lambda _ \varepsilon \Rightarrow $ some $ \psi \in {\mathbf G}^\dagger $  -  forces that for each $ i < \kappa _\varepsilon $, $  \name{g }' _{ \gamma _i (\varepsilon )} \restriction \xi = \nu _i \restriction \xi$\},
                
                \item[$(\beta)$] for $ \bar{ \nu } \in \Xi ^\dagger _ \varepsilon $  let $ \mathbf{G} ^\dagger _ {\bar{ \nu }} $  be like $ \mathbf{G} ^\dagger  $ above,  (it is uniquely determined by   $ \bar{ \nu } $), 
                
                \item[$ (\gamma ) $]  let  $ \Xi^+_\varepsilon = \{ \bar{ \nu } \in \Xi^\dagger_\varepsilon:$  there is a subset $ \mathbf{G} ^+$ of $\mathbf{j} _{ \varepsilon }(\mathbb{P} _{1, \beta (*)}) $  extending $\mathbf{G}^\dagger _ {\bar{\nu}}$  and is  generic over $\bbA_{\varp} =\bfj_{\varp}''(N_{\varp})$ such that $\bbA_{\varp}[\bfG^{+}] = \cH(\chi_{\varp}) \}$ (so $\bfj_{\varp}^{-1}(\bfG^{+})$ is a subset of $\bbP_{1, \beta(*)} \cap N_{\varp}$ generic over $N_{\varp}$ and the Mostowski collapse of $N_{\varp}[\bfj_{\varp}^{-1}(\bfG^{+})]$ is $\cH(\chi_{\varp})) .$  
            \end{enumerate} 
            
            \item[(b)]
            
            \begin{enumerate} 
                \item[$ (\alpha)$] 
                Let $ \Xi ^  \bullet_ \varepsilon $ be the set of pairs $ ( \bar{ \nu }, \mathbf{G} ^+)= \exnu$   such that: $ \bar{ \nu } \in \Xi ^\pigyon_ \varepsilon $ and $ \mathbf{G} ^ + \subseteq  \mathbf{j} _{ \varepsilon }(\mathbb{P} _{1, \beta (*)})$  extending $\bfG_{\overline{\nu}}^{\dagger}$ and is generic over $\bbA_{\varp}.$
                
                \item[$ (\beta )$]  We may write $ (\bar{ \nu } , \mathbf{G} ^+_{\bar{ \nu } })$  or just $ \bar{ \nu } $  though actually $ \bar{ \nu } $ does not determine $ \mathbf{G} ^+_{ \bar{ \nu} }$; but of course $ \bar{ \nu }$  determines $ \mathbf{G} ^\dagger  _{\bar{ \nu } }=  \mathbf{G} ^+_{\bar{ \nu }} \cap \mathbf{j} _{ \varepsilon }(\mathbb{P} '_ { \gamma (*)})$. 
                    
                \item[$ (\gamma )$] Note that $ \bar{ \eta }^ \varepsilon $ belongs to $ \Xi^\pigyon _ \varepsilon $  and even $ (\bar{ \eta } ^ \varepsilon , \mathbf{G} ) \in \Xi ^+_ \varepsilon$ when $ \mathbf{G} =  \mathbf{G} ^+_\varepsilon  $.  
            \end{enumerate} 
        \end{enumerate}
    \end{enumerate}
    
    [Why?  See \cite[3.28-3.32=Le53-Le67]{Sh:1126} and \cite[4.27 = Le70]{Sh:1126}]. 
    
    By the assumption toward contradiction, $\circledast_0$, and $\bbP'_{\ggk}$ being $(< \lambda)$-strategically complete, recalling $\boxplus_1$, there are $\zeta(*),p^{**}$ and $p^+$  (recall $p^* \in \bbP'_{\ggk} \lessdot \bbP_{1,\beta(*)}$  
    is from $ \circledast_0 $ and $ \mathbb{P} _{0, \beta (*)}$ being  dense  in  $  \mathbb{P} _{1, \beta (*)} $) such that:
    
    \begin{enumerate}
        \item[$\circledast_5$]   
        
        \begin{enumerate}
            \item[(a)]  $p^* \le p^{**} \in \bbP'_{\ggk}$ and 
            $p^+ \in \bbP_{0,\beta(*)}$ satisfies   $\bbP_{1,\beta (*)}  \models ``p^{**} \le p^+"$;   (we may add that $ \mathbb{P} '_{\ggk} \models \lqq p^{**}\le  \phi " \Rightarrow \phi , p^+$  are compatible in $\mathbb{P} _{1, \beta (*)}$),
            
            \item[(b)]  $\zeta(*) < \lambda,$
            
            \item[(c)]  $p^{**} \Vdash_{\bbP'_{\ggk}} 
            ``f^*(\zeta) < \name f(\zeta)$  whenever $\zeta(*) \le \zeta < \lambda"$  where $f^*$ is from $\circledast_2,$
            
            \item[(d)]   if $\gamma \in \text{ Dom}(p^+)$ then $\eta^{p^+(\gamma)}$ is an object (not just a $\bbP_{0,\gamma}$-name) and has length $\ge \zeta(*)$ (recall that $\eta^{p^+(\gamma)}$ is the trunk of the condition $p^+(\gamma)$,  see clause $(\alpha)(b)$ of Definition \ref{z23}(1)).
        \end{enumerate}
    \end{enumerate}
    
    Note that possibly Dom$(p^+) \nsubseteq \cup\{v_\varepsilon:\varepsilon < \lambda\}$.  Choose $\varepsilon (*) $  such that:
    
    \begin{enumerate} 
        \item[$\circledast '_5 $]   $\varepsilon(*) < \lambda$ satisfies  $\lambda_{\varepsilon(*)} > \zeta(*) + |\dom(p^{+} )|$ and $ \dom(p^{+} ) \cap \cup \{ v_ \varepsilon : \varepsilon < \lambda \} \subseteq v_{\varepsilon (*)} $ and $\gamma \in \dom(p^+) \Rightarrow \varepsilon(*) > \ell g(\eta^{p^+(\gamma)})$;  
    \end{enumerate}  
    
    Recall clause (d) of $\circledast_5$ and $|\dom(p^+)| < \lambda$ as $p^+ \in \bbP_{0,\beta(*)}$ and $\bbP_{0,\beta(*)}$ is the limit of a $(< \lambda)$-support iteration.
    
    By $\circledast_2$ we can add $(\exists \zeta)[\lambda^-_{\varepsilon(*)} \le \zeta < \lambda_{\varepsilon(*)} < f^*(\zeta)]$.  Our intention is to find $q \in \bbP_{0,\beta(*)}$ above $p^+$ which (in $\bbP_{1,\beta(*)}$)  is above some $q' \in \bbP'_{\ggk}$ which is  $( N_ {\varepsilon (*)}, \mathbb{P} '_\ggk ) $-generic, that is it  forces $\name{\mathbf G}_{\bbP'_{\ggk}}$ to
    include a generic subset of $(\bbP'_{\ggk})^{N_{\varepsilon(*)}}$ hence is induced by some $\bar\nu$ as in $\circledast'_4$, recalling $\circledast_4(b)$.  Toward this in $\circledast_6$ below the intention is that 
    $p^+_{\kappa  (\varepsilon (*))}$ will serve as $q$.
    
    Let $\exi(*) = \exi_{\varepsilon(*)}$ and $\gamma_i$ for $i <\exi(*)$ be such that\footnote{This is used in $\boxplus_3$ and the proof of $(*)_6$.  Not to be confused with $\bar\gamma^\varp$ of $\circledast_3(c)$.}  $\langle \gamma_i:i < \kappa  (*)\rangle$ list $\{\beta^*_i:i \in  v_{\varepsilon(*)}\} =   \cU_{**} = N_{\varepsilon ( * )}  \cap {\cU }   _*  $ in increasing order; recall
    $\cU_* = \{\beta^*_i:i < \ggk\}$ and $i<j < \ggk \Rightarrow \beta^*_i < \beta^*_j$ and $v_{\varepsilon(*)} \subseteq \ggk$ has order type $\exi(\varepsilon(*))$   so $ \gamma _i = \gamma _i(\varepsilon (*))$  from $ \circledast _3$. Next let $\gamma_{\kappa (*)} =  \exi(*)$ so $\{\mathbf
    j_{\varepsilon(*)}(\gamma):\gamma \in v_{\varepsilon(*)}\} = \kappa (*)
    = \mathbf j_{\varepsilon(*)}(\ggk)$.  Recall that $\ggk =  \cf(\kappa) > \lambda,\otp(v_{\varepsilon(*)}) = \otp(N_{\varepsilon(*)} \cap \ggk)$ hence $N_{\varepsilon(*)} \models `` \ggk (*)=  \kappa _{\varepsilon ( *)}$  is a regular cardinal $> \lambda_{\varepsilon(*)}$"; hence:     
    
    \begin{enumerate} 
        \item[(*)] $\exi(*)$ is 
        really a regular cardinal so call it $\sigma$.  
    \end{enumerate} 
    
    Now we define a game $\Game$ as follows\footnote{The idea is to scatter the $\eta^{\varepsilon(*)}_{\gamma_i}$'s.  Why not use the original places? as then we shall have a problem in $\circledast_6$; the scattering is helpful because we are relying on \ref{z35} and \ref{z38}.}:
    
    \begin{enumerate} 
    \item[$\boxplus_2$]   
    
    \begin{enumerate}  
        \item[(A)]  each play lasts $\exi(*)+1 = \sigma +1$ moves and in the $i$-th move: 
        
            \begin{enumerate}  
                \item[(a)]   if $i=j+1$ the antagonist player chooses $\xi_j = \xi(j) < \sigma$ such that $j_1 < j \Rightarrow \zeta (j_1) < \xi(j)$, 
                
                \item[(b)]  then, if $i=j+1$ the 
                protagonist chooses $\zeta_j = \zeta(j) \in (\xi(j),\sigma)$, \underline{but} there are more restrictions implicit in $\boxplus_3$ below, 
                
                \item[(c)]  in any case  (that is, also in the cases $ i \le  \sigma $ is a limit ordinal or zero) the protagonist also chooses $p^+_i,\bar\nu^i$ such that $\boxplus_3$ below holds.  
            \end{enumerate}   
            
            \item[(B)]   in the end of the play the protagonist wins the play \underline{iff} he always has a legal move and in the end:     
            
            \begin{enumerate}   
                \item[(a)]  $ p^+_ \sigma  $ is $ ( {\mathbb   P} '_{\gamma ( * ) }, N_{\varepsilon (*)} ) $-generic,  note the condition is not a member of the same forcing, so we mean that $ p^+_ \sigma$ forces (for $ \mathbb{P} _{1, \beta (*)})$  that the intersection of the generic with $ \mathbb{P} '_{} \cap N_{\varepsilon (*) }$ is generic over $ N_{\varepsilon (*)}$,  
                    
                \item[(b)] $\{\zeta(i):i <  \sigma \}  \in \bbA_{\varepsilon(*)}$; note that trivially it belongs to $M_{\varepsilon(*)} = \bbA_{\varepsilon(*)}[\mathbf{G}^+ _ \varepsilon ] {\mathscr H} ( \chi _ \varepsilon )$,  see $\circledast_4(c)$. 
                    
                \item[(c)] note that we do not demand that  $\bar{ \nu } ' = \langle \nu _{\gamma  _i}: i < \sigma \rangle $  belongs to $ \Xi ^+_ \varepsilon $, we demand only  that  it belongs to  $ \Xi ^\dagger_ \varepsilon $ and even $\Xi_{\varp}'.$
            \end{enumerate} 
        \end{enumerate} 
    \end{enumerate}   
    
    \underline{where},
    
    \begin{enumerate} 
        \item[$\boxplus_3$]
        
        \begin{enumerate}  
            \item[(a)]  $p^+_i \in \bbP_{0,\gamma_i},$
            
            \item[(b)]  if $j<i$ then $\bbP_{0,\gamma_i} \models ``p^+_j \le p^+_i"$,  
            
            \item[(c)]  if $\gamma \in \cup\{\text{Dom}(p^+_j):j<i\}$ then  $p^+_i \restriction \gamma\Vdash_{\bbP_{0,\gamma_i}} ``\name \eta^{p^+_i(\gamma)}$ has length $ \ge i$  and $\ge \lambda_{\varepsilon(*)}"$ moreover $\name \eta^{p^+_i(\gamma)}$ is an object, $\eta^{p^+_i(\gamma)}$,  
            
            \item[(d)]  $\bbP_{0,\gamma_i} \models ``p^{+} \restriction \gamma_i \le p^+_i"$, 
            ($ p^+ $ is from $ \circledast _5$(a)),    
            
            \item[(e)]  
            
            \begin{enumerate}
                \item[$(\alpha)$] $\bar \nu^i = \langle         \nu_{\gamma_j}:j <
                i\rangle$ and $\nu_{\gamma_j} \in 
                \prod\limits_{\iota < \lambda_{\varepsilon(*)}} \theta_\iota$,
            
                \item[$(\beta)$] $\bfG_{\varp(*), i}^{\dagger \dagger}$ is a subset of $\bbP_{\gamma_{i}} \cap N_{\varp(*)}$ generic over $N_{\varp(*)},$
                
                \item[$(\gamma)$] $\bfG_{\varp(\ast), i}^{\dagger} = \bfj_{\varp(*)}''(\bfG_{\varp(*), i}^{\dagger \dagger})$  so is a subset of $\bfj_{\varp}''(\bbP_{\gamma_{2}}')$ generic over $\bbA_{\varp(*)},$
                
                \item[$(\delta)$] $\nu_{j} = \name{\eta}_{\gamma_{j}}[\bfG_{\varp(*), i}^{\dagger \dagger}]$ for $j < i.$ 
            \end{enumerate}
            
            \item[(f)] for $j<i$ we have  $\nu_{\gamma_j} \trianglelefteq \eta^{p^+_i(\gamma_j)}$ so $p^+_i \restriction \gamma_j \Vdash ``\nu_{\gamma_j} \triangleleft \name g'_{\gamma_j}"$ recalling $\boxplus_1$,  
            
            \item[(g)]   for $j<i$  (recall $\bar\eta^\varepsilon$ is from   $\circledast_3(d)$) we have $ ( \alpha ) $  or $ ( \beta ),$ where: 
                
            \begin{enumerate} 
                \item[$(\alpha)$]  $\nu_{\gamma_j} = 
                \eta^{\varepsilon(*)}_{\gamma_{\zeta(j)}}$ recalling
                $\eta^{\varepsilon(*)}_{\gamma_j}$ is from $\circledast_3(d)$,   
                
                \item[$(\beta)$]  $\gamma_j \in \text{ Dom}(p^{+})$ and $\{\iota < \lambda_{\varepsilon(*)}:\eta^{\varepsilon(*)}_{\zeta(j)}(\iota) \ne \nu_{\gamma_j}(\iota)\}$ is a bounded subset of $\lambda_{\varepsilon(*)}$.
            \end{enumerate}  
            
            \item[(h)]  $ p^+_i $ is an upper bound of $\bfG_{\varp(*), i}^{\dagger \dagger}$ hence is a $ (N_ {\varepsilon(\ast)}, \mathbb{P}_{\gamma_{i}}')$-generic  in the natural sense, (actually follows from clause (g), see later in $ \odot $ in the beginning of the proof  of \ref{a35}).
        \end{enumerate} 
    \end{enumerate}  
    
    We shall prove, 
    
    \begin{enumerate}
        \item[$\circledast_6$]  in the game $\Game$:  
        
        \begin{enumerate}
            \item[(a)] the antagonist has no winning strategy,  
            
            \item[(b)]  at  stage $ i $, if $ \langle \zeta ( j ) : j  < i \rangle \in  \mathbb{A} _ \varepsilon $  then the protagonist has a legal move, moreover for any $\zeta(i) \in (\xi(i),\sigma)$ large enough the protagonist can choose it.
        \end{enumerate}
    \end{enumerate}
    
    \underline{Why $\circledast_6$ suffice}?
    
    By clause (a) of $\circledast_6$ we can choose a play $\langle(\xi(i),\zeta(i),p^+_i,\bar\nu^i, \bfG^{\dagger \dagger}_{\varp(*), \sigma}, \bfG_{\varp(*), \nu}^{\dagger}) :i \le \sigma\rangle$ in which the protagonist wins.  Recalling
    $\bbP'_{\ggk} \lessdot \bbP_{1,\beta(*)}$ and
    $\bbP_{0,\beta(*)}$ is a dense sub-forcing of $\bbP_{1,\beta(*)}$, clearly,
    
    \begin{enumerate}
        \item[$\circledast_7$]  there is $p$ such that:  
        
        \begin{enumerate}
            \item[(a)]  $p \in \bbP'_{\ggk}$, 
            
            \item[(b)]  if $\bbP'_{\ggk} \models ``p \le p'"$  hence $p' \in \bbP'_{ \gamma (*)}$ \then \, $p',p^+_\sigma$ are  compatible in $\bbP_{1,\beta(*)}$,  
            
            \item[(c)]  $p$ is above $p^{**}$ and it forces  that $\name g'_{\gamma_i} \rest \lambda_{\varepsilon(*)} = \nu_{\gamma_ {\zeta (i) }}$ for $i < \sigma $  and $\mathbf{j} _{\varepsilon (*)}( \name{{\mathbf G}} _{ {\mathbb   P}' _{ \beta ( * ) } } \cap N _ {\varepsilon  (*)})  = {\mathbf G}^  \pigyon  _{\langle \nu _{\gamma _i }: i < \sigma  \rangle}  \in \Xi ^ \pigyon_{\varepsilon (*) }$.  
        \end{enumerate}
    \end{enumerate}
    
    Then on the one hand,
    
    \begin{enumerate}
        \item[$\circledast_{7.1}$]  $p \in \bbP'_{\ggk}$ being above $p^{**}$ forces $f^* \restriction [\zeta(*),\lambda) < \name f \restriction
        [\zeta(*),\lambda)$ hence $f^* \restriction
        [\zeta(*),\lambda_{\varepsilon(*)}) < \name f
        \restriction [\zeta(*),\lambda_{\varepsilon(*)})$ recalling  that $\zeta(*) < \lambda_{\varepsilon(*)}$, see $\circledast_5$ and the choice of $\varp(*)$ immediately after $\circledast_5$.  
    \end{enumerate}
    
    On the other hand,
    
    \begin{enumerate}
        \item[$\circledast_{7.2}$]  $\bfG_{\varp(*), \sigma}^{\dagger \dagger}$ is a subset of $\bbP' \cap N_{\varp}$ generic over $N_{\varp}.$
    \end{enumerate}
    
    [Why?  By $\circledast_{2}$(e)($\beta$) and the choice of the play.]
    
    \begin{enumerate}
        \item[$\circledast_{7.3}$]  $p_{\sigma}^{+}$ is an upper bound of $\bfG_{\varp(*), \nu}^{\dagger \dagger}.$
    \end{enumerate}
    
    [Why? By $\circledast_{2}$(h) and the choice of the play.]
    
    \begin{enumerate}
        \item[$\circledast_{7.4}$] $p$ is an upper bound of $\bfG_{\varp(*), \nu}^{\dagger \dagger}$ in $\bbP.$ 
    \end{enumerate}
    
    [Why? By $\circledast_{7}$(b), $\circledast_{7.1}$ and $\circledast_{7.2}$].
    
    \begin{enumerate}
        \item[$\circledast_{7.5}$] $p$ is $(N_{\varp(*), \sigma}, \bbP')$-generic. 
    \end{enumerate}
    
    [Why? If $p_{i}^{+}$ is not an upper bound of $\bfG_{\varp(*), i}^{\dagger \dagger}$ (in $\bbP_{1, \gamma_{i}}$) \underline{then} there are $p_{i}''$ and $r \in G_{\varp(*), i}^{\dagger \dagger}$ such that $\bbP_{1, \gamma_{i}} \models$``$p_{i}^{+} \leq p'$'' and $p', r$ are compatible in $\bbP_{1, \gamma_{i}}.$ As $\bbP_{0, \gamma_{i}}$ is a dense subset of $\bbP_{1, \gamma_{i}}$ and $\bbP_{\gamma_{i}}' \cap N_{\varp(*)}$ is a subset of $\bbP_{1, \gamma_{i}}$ of cardinality $\leq \Vert N_{\varp(*)} \Vert < \lambda_{\varp},$ there is $p'' \in \bbP_{0, \gamma_{i}}$ above $p_{i}'$ which decide the value of $\bbP_{\gamma_{i}}' \cap N_{\varp(*)}$].
      
    As $\name f \in N_{\varepsilon(*)}$ it follows from $\circledast_{7.5}$ that:
    
    \begin{enumerate}
        \item[$\circledast_{7.6}$] $p \Vdash  ``\name f \restriction \lambda_{\varepsilon(*)}$ is a function from $\lambda_{\varepsilon(*)}$ to $\lambda_{\varepsilon(*)}"$.
    \end{enumerate}
    
    Together $\circledast_{7.1}$ and $\circledast_{7.6}$ give a contradiction  by the choice of $f^*$ in $\circledast_2$ and of $\varepsilon(*)$ above  which implies that $ \name{ f }( \zeta ) > f^*(\zeta ) >  \lambda _{\varepsilon (*)}$ for some $ \zeta < \lambda _{\varepsilon (*)}$  hence $\circledast_6$ is enough.  In Lemma \ref{a35} below we show that $\circledast_6$ is true; so we are done.
\end{PROOF} 

\begin{lemma}\label{a35}
    The statement $\circledast_6$ is true.
\end{lemma}

\begin{PROOF}{\ref{a35}} 
    Note that:
    
    \begin{enumerate} 
        \item[$\odot$] in $ \boxplus _3,$  clause (h) follows. 
    \end{enumerate} 
    
    [Why?  By \cite[3.93 = Le70]{Sh:1126}, particularly part (5), we have $\circledast_{4}$(c) and the choice of the $\bfG_{\varp}^{\dagger \dagger}$ after $\circledast_{2}.$ In particular recall that  \cite[3.43 = Le70]{Sh:1126} says: 
    
    \begin{itemize}
        \item[$\boxplus$] If $\bfm$ is reasonable (see \cite[2.13 = Le36(3)]{Sh:1126}) \underline{then} for every $p \in \bbP_{\bfm}$ and $s \in \dom(p) \cap M_{\bfm},$ for every large enough $t \in M_{\bfm}$ we have $p \Vdash_{\bbP_{\bfm}}$``$\name{f}_{p(s)} \leq \name{\eta}_{t} \mod J_{\lambda}^{\rm{bd}}$''].
    \end{itemize} 

    Let us prove $\circledast_6$; first, assuming clause (b)  which is proved below,  for clause (a)  choose any strategy {\bf st} for the antagonist and fix a partial strategy {\bf st}$'$ for the protagonist choosing $(p^+_i,\bar\nu^i)$ depending on the previous choices and $\xi (i) < \exi_{\varp(*)}$ such that it is a legal move if relevant and possible.  So the only freedom left for the protagonist is to choose the $\zeta(i)$.  So (recalling $\boxplus_2(A)(a)$) we have in $\mathbf V$ a function $F:{}^{\sigma >}\sigma \rightarrow \sigma$ (so $F$ depends on  $ {\bf st} $ and  ${\bf st}' $)  such that:  
    
    \begin{enumerate}
        \item[$(*)_F$] playing the game such that the antagonist uses {\bf st} and the protagonist uses {\bf st}$'$, arriving at  the $i$-th move,
        $\bar \zeta = \langle \zeta(j):j < i\rangle$ is well defined and if $\bar{ \zeta}   \in N _ {\varepsilon (*)}$  \underline{then}  for the protagonist any choice $\zeta_i \in (F(\bar \zeta),\sigma) \cap
        \cU_{**}$ is legal.
    \end{enumerate}
    
    Note that $F$  belongs to $ {\mathscr H} ( \chi _ \varepsilon)$ unlike $ p^+_ \varepsilon, \bar{ \nu } ^ \varepsilon $.  Now we have to find an increasing sequence $\bar\zeta =  \langle \zeta(i):i < \sigma \rangle$ from $\bbA_{\varepsilon(*)}$ not
    just from $M_{\varp(*)} = \cH(\chi_{\varp(*)})^{\bfV}$ such 
    that $F(\bar \zeta \rest i) < \zeta(i) < \sigma$ and
    $\bar\zeta \in \bbA_{\varepsilon(*)}$. Why possible? As $F \in \cH(\chi_{\varepsilon(*)})$
    and $\cH(\chi_{\varepsilon(*)}) = M_{\varepsilon(*)} = \bbA_{\varepsilon(*)}[\mathbf G^+_{\varepsilon(*)}]$ where $\mathbf G ^+ _{\varepsilon(*)}$ is a subset  of $\mathbf j_{\varepsilon(*)}(\bbP_ {1, \ggk}) \in
    \bbA_{\varepsilon(*)}$ generic over $\bbA_{\varepsilon(*)}$ and $\mathbf j_{\varepsilon(*)}(\bbP_{0,\beta(*)})$ satisfies the 
    $\lambda^+_{\varepsilon(*)}$-c.c. and 
    $\sigma = \text{ cf}(\sigma) > \lambda_{\varepsilon(*)}$ this\footnote{In fact   $ \mathbf{V} \models \lqq \mathbb{P} '_ \kappa $  satisfies the $ \kappa $-c.c."   suffices.}  is possible.  That is, there is a $\mathbf j_{\varepsilon(*)}(\bbP_{0,\beta(*)})$-name
    $\name F_* \in \bbA_{\varepsilon(*)}$ such that $F = \name F_*[\mathbf G ^+ _{\varepsilon(*)}]$ and we define in  $\bbA_{\varepsilon(*)}$ the function  $F':{}^{\sigma >}\sigma \rightarrow \sigma$ by  $F'(\langle \zeta(j):j < i  \rangle ) = \sup\{\xi +1:\xi \in \{\zeta(j)+1:j<i\}$ \underline{or} $\xi < \sigma$ and $\nVdash_{\mathbf j(\bbP_{0,\beta(*)})} ``\name F(\langle \zeta(j):j<i\rangle) \ne 
    \xi" \}$; clearly this is O.K.
    
    We are left with proving $\circledast_6(b)$.
    
    \underline{Case 1}:  $i=0$.
    
    Let $p^+_0 = p^+ \restriction \gamma_0$.
    
    \underline{Case 2}:  $i$ limit.
    
    By clauses (b) and (c) of $\boxplus_3$, there is 
    $p^+_i \in \bbP_{0,\gamma_i}$ which is an
    upper bound (even l.u.b.) of $\{p^+_j:j <i\}.$ Note that $\overline{\nu}_{i} = \langle \nu_{j}: j <i \rangle$ and it satisfying $\boxplus_{2}(B)(b)$ and $\boxplus_{3}$(e)($\alpha$)(f)(g), so by $\boxplus_{1}(b)^{+}$ there is $\bfG_{\varp(*), i}^{\dagger}$ as required in $\boxplus_{3}(\beta), (\gamma).$ 
    
    So we are done with Case 2. 
    
    \underline{Case 3}:  $i=j+1$ and $\gamma_j \notin \dom(p^+)$.  
    
    Clearly $\gamma_i$ is in $\cU_*$ the successor of 
    $\gamma_j$ and $(\exists \iota)(\gamma_j = \beta^*_\iota \wedge \iota \in v_{\varp(*)})$. As in case 4 below but easier by the properties of the iteration  and \cite[\S3C]{Sh:1126}.
    
    \underline{Case 4}:  $i=j+1$ and $\gamma_j \in \dom(p^+)$. Again $\gamma_i$ is in $\cU_*$ the successor of $\gamma_j$ and $(\exists \iota)(\gamma_j = \beta^*_\iota \wedge \iota \in v_{\varp(*)})$. 
    
    First we find $p'_j$ such that:
    
    \begin{enumerate}
        \item[$\circledast_8$]  
        
        \begin{enumerate}
            \item[(a)]  $p^+_j \le p'_j \in \bbP_{0,\gamma_j}$, 
            
            \item[(b)]   if $\gamma \in \dom(p^+_j)$ then
            $p'_j \restriction \gamma \Vdash ``\ell g(\name\eta^{p'_j(\gamma)}) > i(*)= \sigma "$ (see $\boxplus_3(c)$),   
            
            \item[(c)]  $p'_j$ forces \footnote{recall 
            that $\eta^{p^+(\gamma_j)}$ is an object, not a name and $p^+_j$ is $(N_{\varepsilon(*)},\bbP_{\gamma_j})$-generic} a value to the pair $(\eta^{p^+(\gamma_j)},\name f^{p^+(\gamma_j)}  \restriction \lambda_{\varepsilon(*)})$; we call this pair $q_j = (\eta^{q_j},f^{q_j})$.
        \end{enumerate}
    \end{enumerate}
    
    [Why? This should be clear.]
    
    Second,
    
    \begin{enumerate}
        \item[$\circledast_9$]  $p^{+} _j$ hence $p'_j$ is $(N_{\varepsilon(*)},\bbP'_{\gamma_j})$-generic and $\langle \nu_{\gamma_{j(1)}}:j(1) < j\rangle$ induces the generic.
    \end{enumerate}
    
    [Why? By clause (h) of $  \boxplus _2 $, see $ \odot $ above.   Alternatively  As in the proof of $\circledast''_7$ of Lemma \ref{a32} when we assume that we have carried the induction, by $\boxplus_2$, clause (g) and $\circledast_4$].
    
    Now,
    
    \begin{enumerate}
        \item[$\circledast_{10}$] 
        
        \begin{enumerate}
            \item[(a)]  $f^{q_j} \in (\prod_{\zeta < \lambda_{\varp(*)}} \theta_\zeta)^{\bbA_{\varepsilon(*)}[ \mathbf{G} ^+_{\varepsilon (*)} ]}$; recalling that $ f^{q_j}$  is from clause (c) of $\circledast_8$. 
            
            \item[(b)]   for every large enough $\zeta \in (\xi(i),\sigma)$ we have:
            
            \begin{enumerate}
                \item[$\bullet$]  $f^{q_j} \le \eta^{\varepsilon(*)}_\zeta \mod J^{\bd}_{\lambda_\varepsilon}$. 
            \end{enumerate}
        \end{enumerate}
    \end{enumerate}
    
    [Why?  Clause $\circledast_{10}$(a) holds because $f^{q_j} \in (\prod\limits_{\zeta < \lambda_{\varp(*)}} \theta_\zeta)^{\mathbf V}$, hence belongs to $\cH(\chi_{\varp(*)})$ which is the universe of $M_{\varp(*)}$ so $f^{q_j} \in M_{\varp(*)}$.  But $M_{\varp(*)} = \bbA_{\varp(*)}[\mathbf{G}^+  _ {\varepsilon (*)}]= {\mathscr H} (\chi _ {\varepsilon (*)}, \in )$ and $ \bar{ \eta }^{\varepsilon (*)}= \langle \mathbf{j} _{\varepsilon(*)} (\name{ \eta } _ \gamma ): \gamma \in \cU_{\ast} \cap N_{\varepsilon (*)}\rangle$; recalling  $\bar\eta^{\varp(*)}$ is a generic for $\mathbf j_\varp(\bbP'_{\ggk})$. 
    
    For clause $\circledast_{10}$(b) recall $(*)_4(b)$.  Hence $ N_{\varepsilon (*)} $ satisfies the parallel statement,  so $N_{\varepsilon(*)}$ satisfies: if we force by $ \mathbb{P} _ {}$ then $ \{ \name{ \eta } _ \gamma : \gamma \in  {\mathscr U } _* \cap N_{\varepsilon(*)} \} $  is cofinal in $ (\Pi _{\varepsilon < \lambda_{\varepsilon(*)}}  \theta _ \varepsilon, \le _{J^{\bd }_{\lambda _{\varepsilon(*)}}})$. Note that for every $\varp, \otp(\cU_{\ast} \cap N_{\varp})$ has cofinality $> \lambda_{\varp}$ by the choice of $N_{\varp}.$ 
        
    This  is a crucial point: this is justified by clause  (A)(e) of \ref{z38}.   
    
    Applying $ \mathbf{j} _{\varepsilon(*) }$ and recalling $ \mathbb{A} _{\varepsilon (*)}[\mathbf{G}^+ _{\varepsilon (*)}] = {\mathscr H} ( \chi _{\varepsilon (*) })$  we are done proving $ (*)_{10}$]. 
    
    Now we choose $\zeta(j)  > \sup \{ \zeta ( j_1 ) : j_1 < j\} $ as in clause (b) of $\circledast_{10}$ and $\nu_j  = \eta _{\zeta (j)}$;  so here we obey the promise \lqq for every large enough $ \zeta(i)$".  
    Next choose $p^+_i \in \bbP'_{\ggk}$ such that $p^+_i \rest \gamma_j = p'_j,\eta^{p^+_i(\gamma_i)} = \nu_j$ and $f^{p^+_i  (\gamma_j )} \rest [\lambda_\varepsilon,\lambda) = f^{p^+(\gamma_ j)} \rest [\lambda_\varepsilon,\lambda)$ and $ \nu _{\zeta (j)} \triangleleft f^{p^+_i( \gamma _j )}$". 
    
    Lastly, we choose $p^+_i $ above $p'_i $ as in 
    the proof of Case 2, so we have finished Case 4.   
    
    We have carried the induction hence proved $\circledast_6(b)$ so we are done proving \ref{a35}.
\end{PROOF}

\begin{discussion}\label{a41}\ 

    (1) The reader may justly wonder why we use 
    $\mathbb{A}' = \mathbf \mathbb{A} [\name{\bar g}'] = \mathbb{A} [\name{\bar g} \rest \cU_*]$ rather than simply $\mathbb{A} [\name{\bar g}]$.  
    Of course, nothing is lost by it, but why the extra complication?
    
    (2) The answer is that we are committed to $ p^+ $  so $\mathbb{P} _{0, \beta (*) }$, and it is not clear why it can be increased to a condition $ q $  which is $ (N_ \varepsilon (*))$-generic and forces the desired  statement (i.e. contradicting  \lqq $ f $ dominates  $ \langle f^* _  \alpha  : \sigma < \mu \rangle $"). We succeed  to do this using $ \bar{ \nu } ' $ which is almost equal  to a suitable sub-sequence of $ \eta ^{\varepsilon(*)}.$ So  during the proof we  used: if $\zeta(i) \in \cU_*$ is increasing with $i < \ggk$ then also $\langle \name  \ttg '_{\zeta(i)}:i < \kappa\rangle$ is generic over $\mathbf V$ for the sub-forcing of $\bbP_{1,\beta(*)}$ generated by $\name{\bar g} \rest \cU_*$; see $\circledast''_7$ inside the proof of $\circledast_6$ inside \ref{a35}. But using $\cU_* = \beta(*)$, we do not know this.
    
    (3) Now in the parallel case for $\lambda = \aleph_0$ with FS-iteration with full memory, such claim is true, see \S0.  
    
    (4) But we do not know the parallel of (3) for $\lambda$, so we use a substitute using $\cU_*$, i.e. $\bbP'_\kappa$.
\end{discussion}

\begin{claim}\label{a44}
    In \ref{a32} we can add $\bbP_{\kappa}$``$\lambda$ is supercompact''.
\end{claim}

\begin{PROOF}{\ref{a44}}
    Recall our forcing actually is $\bbR \times \name{\bbP}^{1} \ast \name{\bbP}^{2},$ where $\bbR$ is the preparatory forcing build from Laver's diamond, so $\bbR = \bbR_{\lambda},$ limit of the Easton support iteration $\langle \bbR_{\alpha}^{1}, \name{\bbR}_{\beta}^{0}: \alpha \leq \lambda, \beta < \lambda \rangle, \bbP_{1}$ forces $\mathfrak{d}_{\lambda} = \mathfrak{b}_{\lambda} = \mu,$ and $\name{\bbP}^{2}$ is $\bbP'$  use in \ref{a32} (all over $\bfV_{0}$). 
    
    For every $\mu$ we can find (in $\bfV_{0}$) a transitive class $\bfM, \bfM^{< \chi} \subseteq M$ and elementary embedding $\bfj: \bfV \to \bfM,$ with critical cardinal $\lambda.$ 
    
    Repeating the proof of preservation of supercompactness it is enough find an upper bound for $\bfj''(\bfG_{\bbP}^{2}),$ what is done as in the proof \ref{a35}, with $\gamma_{2} = i.$   We elaborate in \cite[4.28 = Le23]{Sh:1126} 
\end{PROOF}

\bibliographystyle{amsalpha}

\bibliography{shlhetal}

\end{document}